\documentclass[10pt]{article}
\usepackage{amssymb}
\usepackage{latexsym}
\usepackage{euscript}

\newcommand{\pp}{\mathbb{P}}
\newcommand{\qq}{\mathbb{Q}}
\newcommand{\cc}{\mathbb{C}}
\newcommand{\rr}{\mathbb{R}}

\newcommand{\zz}{\mathbb{Z}}
\newcommand{\sss}{\mathbb{S}}
\newcommand{\rot}{\mathbb{T}}

\newcommand{\nn}{\mathbb{N}}

\DeclareFontFamily{OT1}{rsfs}{}
\DeclareFontShape{OT1}{rsfs}{n}{it}{<-> rsfs10}{}
\DeclareMathAlphabet{\curly}{OT1}{rsfs}{n}{it}

\newcommand{\cdbar}{\overline{\partial}}

\newcommand{\im}{\mathrm{im \,}}

\newcommand{\Sym}{\mathrm{Sym}}

\newcommand{\Crit}{\mathrm{Crit}}

\newcommand{\Pic}{\mathrm{Pic \,}}
\newcommand{\ind}{\mathrm{Index \,}}

\newcommand{\Fibre}{\mathrm{Fibre}}

\newcommand{\Jac}{\mathrm{Jac}}

\newcommand{\mgbar}{\overline{M}_g}

\newcommand{\Hilb}{\mathrm{Hilb}}
\newcommand{\cok}{\mathrm{cok}}

\newcounter{Universal}[section]
\renewcommand{\theUniversal}{\thesection.\arabic{Universal}}

\newenvironment{Plain}{\refstepcounter{Universal} \par \vspace{0.5cm}
\noindent {\bf (\theUniversal)}\ }{\par \vspace{0.5cm}}
\newenvironment{Italic}{\refstepcounter{Universal} \par \vspace{0.5cm}
\noindent {\bf (\theUniversal)}\ \it}{\par \vspace{0.5cm}}  
  
\newenvironment{Thm}{\begin{Italic}{\sc Theorem: }}{\end{Italic}}
\newenvironment{Prop}{\begin{Italic}{\sc Proposition:} }{\end{Italic}}
\newenvironment{Defn}{\begin{Italic}{\sc Definition: } }{\end{Italic}}  
\newenvironment{Cor}{\begin{Italic}{\sc Corollary: }}{\end{Italic}}
\newenvironment{Lem}{\begin{Italic}{\sc Lemma: }}{\end{Italic}}
\newenvironment{Conj}{\begin{Italic}{\sc Conjecture: }}{\end{Italic}}
\newenvironment{Pf}{\par \noindent{\sc Proof:} }{\quad $\blacksquare$ \par
\vspace{0.5cm}}
\newenvironment{Example}{\begin{Plain} \noindent{\sc Example: }}
{\quad $\square$ \end{Plain}} 
\newenvironment{Rmk}{\begin{Plain} \noindent{\sc Remark: }}{\quad
    $\square$ \end{Plain}}

\newcounter{enum}

\newenvironment{Eqn}{\refstepcounter{Universal} $$} {\eqno \mathrm{
(\theUniversal)} $$} 

\title{Serre-Taubes duality for pseudoholomorphic curves}
\author{Ivan Smith\thanks{New College, Oxford, OX1 3BN, England:
    smithi@maths.ox.ac.uk}} 
\date{}

\begin{document}
\maketitle
\thispagestyle{empty}

\begin{abstract}
\noindent  According to Taubes, the Gromov invariants of a symplectic
four-manifold $X$ with $b_+ > 1$ satisfy the duality $Gr(\alpha) = \pm
Gr(\kappa - \alpha)$, where $\kappa$ is Poincar\'e dual to the
canonical class.  Extending joint work with Simon Donaldson, we
interpret this result in terms of Serre
duality on the fibres of a Lefschetz pencil on $X$, by proving an
analogous symmetry for invariants counting sections of
associated bundles of symmetric products.  Using similar methods, 
we give a new proof of an existence theorem for symplectic surfaces in
four-manifolds with $b_+ =1$ and $b_1=0$.  This reproves another theorem
due to Taubes: two symplectic homology projective planes with negative
canonical class and equal volume are symplectomorphic. 
\end{abstract}

\section{Introduction}

Many questions concerning the topology of symplectic manifolds can be
formulated in terms of the existence (or otherwise) of appropriate
symplectic submanifolds.  In four dimensions, this viewpoint has been
especially fruitful given rather general existence theorems for
symplectic surfaces coming from Taubes' theory relating Seiberg-Witten
invariants to pseudoholomorphic curves \cite{Taubes}.  Taubes' results
imply in particular a curious and striking symmetry for counting
holomorphic curves in symplectic four-manifolds with $b_+ > 1$:  the Gromov
invariants for the classes $\alpha$ and $\kappa-\alpha$ are equal up
to sign, where
$\alpha, \kappa \in H_2 (X; \zz)$ and $\kappa$ is Poincar\'e dual to the
canonical class of the symplectic structure.  (As we recall in the next
section, this is a consequence of Serre duality when $X$ is K\"ahler,
and can be seen as a symplectic shadow of the Serre duality theorem.)
Given Taubes' amazing
theorem on the equivalence ``SW = Gr'' this is rather trivial: it is
the translation of a routine symmetry in Seiberg-Witten theory arising
from flipping a spin$^c$-structure $L$ to its dual $L^*$.
Nonetheless, one might hope that there was a direct proof of the
symmetry in the context of Gromov invariants and holomorphic curves,
which are in a real sense the more geometric objects of interest in
the symplectic setting.  This paper is designed, amongst other things,
to outline one route to such a geometric interpretation of the duality.

\vspace{0.2cm}

\noindent In \cite{SKDIS}, in joint work with Simon Donaldson, we
explained how to obtain symplectic surfaces in four-manifolds with
rational symplectic form from
Lefschetz pencils.  These arise as sections of associated fibrations
of symmetric products down the fibres of the Lefschetz pencil;  using
this viewpoint, we were able in particular to give a new proof of
Taubes' theorem that for ``most'' symplectic four-manifolds the
class $\kappa$ can be represented by an embedded symplectic surface.
More precisely (a fuller review will be given in the third
section) we defined an invariant - the ``standard surface count'' -
$\mathcal{I}_{(X,f)}(\alpha)$ as
follows.  Be given a Lefschetz pencil $f: X \dashrightarrow \sss^2$ of
genus $g$ curves 
and construct a fibre bundle $F:X_r(f) \rightarrow \sss^2$ with fibres
the $r=(2g-2)$-th symmetric products of fibres of $f$.  There is a
natural injection $\phi$ from the set of sections of $F$ to $H_2 (X; \zz)$.
Then $\mathcal{I}_{(X,f)}(\alpha)$ counted the holomorphic sections of
$F$ in the unique class $\tilde{\alpha}$ with image
$\phi(\tilde{\alpha}) = \alpha$, and was defined to be zero for
$\alpha \not \in \im(\phi)$.  The main theorems of \cite{SKDIS}, for
manifolds with rational symplectic form, were then: 

\begin{itemize}
\item if $\mathcal{I}_{(X,f)}(\alpha) \neq 0$ then $\alpha$
may be represented by an embedded symplectic surface in $X$;
\item for any $X$ with $b_+ (X) > 1+b_1 (X)$ and any Lefschetz pencil
  $f$ on $X$ (of sufficiently high degree) we have
  $\mathcal{I}_{(X,f)}(\kappa) = \pm 1$. 
\end{itemize}

\noindent The main theorems of this sequel paper are, in these terms, the
following:  again $\mathcal{I}$ will denote an invariant which counts
sections of a compactified bundle of $r$-th symmetric products, for
arbitrary possible $r$.

\begin{Thm} \label{Serreduality}
Fix a symplectic four-manifold $X$ and
a Lefschetz pencil $f$ of sufficiently high degree on $X$:
\begin{itemize}
\item if $b_+(X) > 1+b_1(X)$ we have an equality $\mathcal{I}_{(X,f)}(\alpha)
  = \pm \mathcal{I}_{(X,f)}(\kappa - \alpha)$. 
\item if $b_+(X)=1$ and $b_1(X)=0$ and $\alpha\cdot[\omega]>0, \, \alpha^2 >
  K_X\cdot\alpha$ then $\mathcal{I}_{(X,f)}(\alpha) = \pm 1$.
\end{itemize}
\end{Thm}

\noindent The broad strategy is as in \cite{SKDIS}.  First we set up
enough theory to define the relevant Gromov invariants.  Then we use
almost complex structures adapted to diagonal loci in symmetric
products to obtain symplectic surfaces, and lastly we use almost
complex structures adapted to the geometry of the Abel-Jacobi map to
perform explicit computations.  (The third stage in particular is
conceptually as well as technically more involved here.)  These
results should be considered together with the following

\begin{Thm} \label{IequalsGrmodtwo}
Suppose $(X, \omega)$ contains no embedded
symplectic torus of square zero.  Then there is an equality
$\mathcal{I}_{(X,f)} = Gr_X: H_2 (X; \zz) \rightarrow \zz \rightarrow
\zz_2$ of the
mod two reductions of the standard surface count and the Gromov invariant.
\end{Thm}

\noindent The hypothesis on $X$ is satisfied, for instance, if $K_X =
\lambda [\omega_X]$ for any $\lambda \in \rr^*$.  Here the invariant
$Gr_X$ is that introduced by Taubes,
counting embedded but not necessarily connected holomorphic curves in
$X$.  It would be natural to extend the theorem by taking account of
the signs (to lift to an equivalence of integer-valued invariants) and by
including the delicate situation for tori; we leave those extensions
for elsewhere.   Up to certain universal weights in the definition of
the $\mathcal{I}$-invariants when tori are present, we expect:

\begin{Conj} \label{IequalsGr}
Once $f$ is of sufficiently high degree, $\mathcal{I}_{(X,f)}$ is
independent of $f$ and defines a
symplectic invariant of $(X, \omega)$.  Moreover there is an equality
$\mathcal{I}_{(X,f)} = Gr_X: H_2 (X; \zz) \rightarrow \zz$.
\end{Conj}

\noindent It seems likely, given the stabilisation procedure for
Lefschetz pencils, that the first half of this can be proven directly
without identifying $\mathcal{I}$ and $Gr_X$.  We make a few remarks
on this at the end of the fifth section.   Of
course the conjecture, combined with Taubes' results,
would imply the first theorem.  Rather, the point is that one
should be able to prove the conjecture independently of Taubes'
results, and then (\ref{Serreduality}, part (i)) yields a new proof of the
symmetry of the Gromov invariants.  In this paper, we shall concern
ourselves with the properties of the $\mathcal{I}$-invariants, but let us point
out that along with (\ref{IequalsGrmodtwo}) we can obtain
holomorphic curves in any $X$:

\begin{Thm}
Let $X$ be a symplectic four-manifold and fix a taming almost complex
structure $J$ on $X$.
If $\mathcal{I}_{(X,f)}(\alpha) \neq 0$ for all sufficiently high
degree pencils $f$ on $X$ then $X$ contains $J$-holomorphic curves in the
class $\alpha$.
\end{Thm}

\noindent  As a technical remark, note that this enables one to find
symplectic surfaces in $X$ even if the symplectic form is not
rational.  One might hope that the equivalence
$\mathcal{I} = Gr$, together with the conjectural
routes to establishing $SW = \mathcal{I}$ as proposed by
Dietmar Salamon, would give a more intuitive framework for  
Taubes identification $SW = Gr$.
In our picture, the symmetry (\ref{Serreduality}, (i)) will
arise entirely 
naturally from Serre duality on the \emph{fibres} of the Lefschetz fibration;
the key geometric ingredient is the fact that the $(2g-2-r)$-th
symmetric product of a Riemann surface fibres over the Jacobian with
exceptional fibres precisely over an image of the $r$-th symmetric
product, when $r < g-1$.  For technical reasons, we will make use of a
strengthening of
this observation which gives us enough control to work in families
(\ref{nicefibres}). 
This stronger result will follow easily from results, due to Eisenbud
and Harris, in the
Brill-Noether theory of Riemann surfaces;
one appeal of the current proof is that such results become of
relevance in four-dimensional symplectic geometry. 

\vspace{0.2cm}

\noindent The symmetry of the Gromov invariants is false for
symplectic manifolds with $b_+=1$.  Indeed, if $X$ is minimal and has
$b_+=1$, and if in addition $Gr(\alpha)$ and $Gr(\kappa-\alpha)$ are
both non-trivial, then necessarily $K_X^2=0$ and $K_X = n\alpha$ for
some integer $n$.  This is an easy consequence of properties of the
intersection form on such four-manifolds; for K\"ahler surfaces, if
both $Gr(\alpha)$ and $Gr(\kappa-\alpha)$ are non-trivial then the
holomorphic curves in the two homology classes come from sections of
line bundles which may be tensored to produce a non-trivial element of
$H^0(K_X)$, forcing $b_+ > 1$.  Taubes in fact proves the theorem
under the weakest 
possible constraint $b_+>1$.  We shall assume throughout the bulk of
the paper that $b_+>1+b_1$; this weaker assumption simplifies the
arguments, and keeps the geometry to the fore.  At the end of the
paper we shall sketch how to improve the arguments to hold in case
$b_+>2$.  This still falls short of Taubes, and superficially at least
the Hard Lefschetz theorem plays a complicating role.  Presumably
sufficient ingenuity would cover the missing case $b_+=2$, but the
author could not find an argument.
In any case, rather than being sidetracked we hope to emphasise the
key geometric ingredients of the new proof.  After all, the theorem
itself already has one beautiful and very detailed exposition, thanks
to Taubes, and our intention is to supplement and not supplant the
gauge theory\footnote{``Il n'y a point de secte en g\'eom\'etrie'',
  Voltaire (Dictionnaire Philosophique)}. 
 
\vspace{0.2cm}

\noindent The non-vanishing result (\ref{Serreduality}, part (ii)) -
which is described at the end of the paper - 
implies in particular an existence theorem for symplectic surfaces in
four-manifolds with $b_+ = 1$.
Such results are well-known, and go back to McDuff
\cite{McDuff:DeftoIso}; similar work has been done by T.J.Li and
A.K.Liu \cite{Li-Liu}.  Each of these earlier proofs has relied on
wall-crossing formulae for Seiberg-Witten invariants;  our
arguments are ``more symplectic'' and may
cast a new light on the relevant geometry.  We remark that here it is
important to
use a definition of $\mathcal{I}(\alpha)$ in which we
cut down a positive dimensional moduli space by intersecting the image
of an evaluation map with appropriate divisors (corresponding to
forcing holomorphic curves to pass through points in the four-manifold
$X$); we will explain this more properly below.

\vspace{0.2cm}

\noindent The proofs of both parts of (\ref{Serreduality}) run along
similar lines 
to the proof of the main theorems of \cite{SKDIS}, and much of the
technical material is already present in that paper.  As before,
monotonicity of the fibres of $X_r(f)$ enables us to use elementary
machinery from the theory of pseudoholomorphic curves, so the proofs
are not too hard.  Using the
results of the two papers together, it now becomes
possible to re-derive some of the standard structure theorems for
symplectic four-manifolds from the perspective of the existence of
Lefschetz pencils.  Here are some sample results,
whose proofs are well-known: for completeness we briefly recall the
arguments in the last section of the paper.

\begin{Cor}\label{applications}
Let $X$ be a symplectic four-manifold with $b_1 = 0$. 
\begin{itemize}
\item (Taubes) If $X$ is minimal with $b_+ > 1$ then $2e(X)+3\sigma(X) \geq
  0$.  In particular, manifolds such as $K3 \sharp K3 \sharp K3$ admit
  no symplectic structure.
\item (Taubes) A homology symplectic projective plane with $K_X \cdot
  [\omega] < 0$ 
  is symplectomorphic to $(\cc \pp^2, \mu\omega_{FS})$ for some $\mu >0$.
\item (Ohta-Ono) More generally, if $c_1 (X) = \lambda [\omega]$ for
  some $\lambda \in \rr_{>0}$ 
  then $X$ is diffeomorphic to a del Pezzo surface.  
\item (Li-Liu) If $X$ is minimal with $b_+ = 1$ and $K_X^2 > 0, \, K_X \cdot
  \omega > 0$ then the canonical class contains symplectic forms.
\end{itemize}
\end{Cor}

\noindent  The examples, though not exhaustive, serve also to
highlight some of the 
profound successes of the gauge theory which remain mysterious from
the perspective of symplectic linear systems:  one such is the role of
\emph{positive scalar curvature} as an obstruction to the existence of
holomorphic curves. 

\vspace{0.2cm}

\noindent Let us remark on three further directions suggested by
\cite{SKDIS} and this paper.  The first concerns non-symplectic
four-manifolds.  Work of Presas \cite{Presas} suggests that symplectic
manifolds with contact boundary should also admit pencils of sections, and
one could hope to complement Taubes' theorems on the Seiberg-Witten
equations on manifolds with 
self-dual forms \cite{Taubes:selfdual} with \emph{existence} statements for
holomorphic curves with boundary.  A second concerns higher dimensional
symplectic manifolds.  For complex three-folds fibred smoothly over
curves, the relative Hilbert schemes are smooth and one can approach
the Gromov invariants of the three-folds through sections of the
associated bundles.  Counting such sections gives invariants of loops
of symplectomorphisms for complex surfaces which refine those of
\cite{Seidel2}.  More
generally, after finitely many blow-ups any symplectic  six-manifold $Z$
admits a map to the projective plane;
a complex curve $\Sigma$ in $Z$ projects to a complex curve $C$ in
$\cc \pp^2$, and 
generically at least $\Sigma$ lies inside the total space of a Lefschetz
fibration over $C$.  This suggests an inductive approach to the Gromov
invariants of $Z$, similar in flavour to work of Seidel
\cite{SeidelVM}.  In a third direction, both \cite{SKDIS} and this paper
concern solutions of the Seiberg-Witten equations on the total space
of a Lefschetz fibration.  An analogous story for the instanton
equations is the subject of work in progress by the author and will be
the topic of a sequel paper.

\vspace{0.2cm}

\noindent \textbf{Outline of the paper:}

\begin{enumerate}

\item Remark: although we will not repeat all details of the
  \emph{local} constructions of \cite{SKDIS}, we will give a (more) coherent
  development of the \emph{global} theory that we require, so the
  paper should be accessible in its own right.

\item In the next section, we explain how the symmetry $Gr_X(\alpha) =
  \pm Gr_X(\kappa-\alpha)$ follows from Serre duality if $X$ is a K\"ahler
  surface with $b_1=0$.  (This motivates various later constructions.)

\item In the third section, we recall the basics of Lefschetz
  fibrations, and prove that given $(X,f)$ the relative Hilbert scheme
  $X_r(f)$ provides a smooth symplectic
  compactification of the family of $r$-th symmetric products of
  the fibres.  (We have tried to illuminate the structure of this space.)

\item In the fourth section, we compute the virtual dimensions of
  moduli spaces of 
  sections of $X_r(f)$  and define an invariant $\mathcal{I}$
  which counts sections in a fixed homotopy class (this requires a
  compactness theorem).  We also give a simple ``blow-up'' formula.

\item In the fifth section, using a natural almost
  complex structure on $X_r(f)$ and Gromov compactness, we show that
  if $\mathcal{I}(\alpha) \neq 0$ then the moduli space of holomorphic
  curves representing $\alpha$ is non-empty for any taming almost
  complex structure on $X$.  We also sketch how to obtain the equivalence
  $\mathcal{I} = Gr_X \  (\mathrm{mod} \ 2)$ for manifolds containing
  no symplectic square zero tori.

\item In the sixth section, we prove the main result
  (\ref{Serreduality}); this involves a short detour into
  Brill-Noether theory and some obstruction computations modelled on
  those of \cite{SKDIS}.  We assume, for simplicity, that $b_+ >
  1+b_1$ or $b_+ = 1, b_1 = 0$ in this section.

\item In the final section, we give the proofs of the applications
  listed above and explain how to extend the arguments of
  \cite{SKDIS} (and in principle Section 6) to the case where $b_+ > 2$.
\end{enumerate}

\vspace{0.4cm}

\noindent \textbf{Acknowledgements:}  Many of the ideas here arose
from conversations with Simon Donaldson, whose influence has been
accordingly pervasive.  Thanks also to Denis Auroux, Eyal Markman,
Paul Seidel and Bernd Siebert for helpful remarks.


\section{Digression on algebraic surfaces}

We shall begin (semantic sensibilities regardless!) with
a digression.  If $X$ is K\"ahler then one can
often compute the Gromov invariants of $X$ directly;  we will review
this, and explain how the symmetry $Gr_X (\alpha) = \pm Gr_X
(\kappa-\alpha)$ emerges in this framework.  Suppose for simplicity that $b_1
(X) = 0$.

\vspace{0.2cm}

\noindent  The key point is that a holomorphic
curve, for the integrable complex structure, is exactly a divisor and
as such gives rise to a section of a line bundle.  Moreover,
generically at least, the locus
of holomorphic sections of a given line bundle yielding singular complex curves
will have positive codimension, whilst those yielding curves with
worse than nodal singularities will have complex codimension at
least two.  It follows that the linear system in which the divisor
moves defines a suitable compactification of the space of smooth
holomorphic curves for computing Gromov invariants.  Again for
conceptual clarity, and since we shall not use the results of this
section later on, we will suppose we are always in this situation.  The
desired invariant itself can be computed as the Euler class of 
an obstruction bundle over the moduli space (cf. \cite{Salamon:SW},
Prop. 11.29).

\begin{Prop} \label{quotientalg}
Let $X$ be a K\"ahler surface with $b_+>1$ and $b_1=0$.  
The Gromov invariants of $X$ manifest the symmetry $Gr(\alpha) = \pm
Gr(\kappa-\alpha)$. 
\end{Prop}

\begin{Pf}
Fix a suitable class $\alpha \in H_2 (X;\zz)$ which we suppose to be Poincar\'e
dual to a class $D$.  We will blur $D$ with the unique holomorphic
line bundle $\mathcal{O}(D)$ with first Chern class $D$.  Our
  assumptions imply that 
we have $H^1 (X, \mathcal{O}_X) = 0$; then $\pp = \pp (H^0 (D))$ is
the moduli space of pseudoholomorphic curves in the homology class
$\alpha$ - there are no other line bundles with the same
first Chern class. 
Write $\Delta \subset \pp \times X$ for the universal divisor with 
projections $\pi: \Delta \rightarrow \pp$ and $p: \Delta \rightarrow
X$.  The obstruction bundle
by definition is given by $R^1 \pi_* \mathcal{O}_{\Delta}(\Delta)$.
Suppose first that the virtual dimension is zero; $\alpha^2 =
K_X\cdot\alpha$.  Then the invariant is the Euler class of the
obstruction bundle.
We have an exact sequence

$$0 \rightarrow \mathcal{O}_{\pp \times X} \rightarrow
\mathcal{O}_{\pp \times X} (\Delta) \rightarrow \mathcal{O}_{\pp \times
  X} (\Delta) |_{\Delta} \rightarrow 0$$

\noindent where the last non-zero term is
$\mathcal{O}_{\Delta}(\Delta)$ by definition of notation.  It follows
that

$$\pi_! \mathcal{O}_{\pp \times X} + \pi_! \mathcal{O}_{\Delta}
(\Delta) \ = \ \pi_! \mathcal{O}_{\pp \times X} (\Delta)$$

\noindent in K-theory, and hence taking total Chern classes that 

\begin{Eqn} \label{samechernclass}
c (\pi_! \mathcal{O}_{\Delta} (\Delta) ) \ = \ c( \pi_!
\mathcal{O}_{\pp \times X} (\Delta) ).
\end{Eqn} 

\noindent Now $\pi: \Delta \rightarrow \pp$ has one-dimensional fibres
and hence 

$$\pi_! \mathcal{O}_{\Delta}(\Delta) \ = \ R^0 \pi_*
\mathcal{O}_{\Delta}(\Delta) - R^1 \pi_*
\mathcal{O}_{\Delta}(\Delta).$$

\noindent Note that the first of these is by deformation theory just
the tangent space to the moduli space $\pp$, whilst the latter is the
obstruction space we require.  It follows, also using
(\ref{samechernclass}),  that

$$c (\mathrm{Obs}) \ = \ c (\pi_! \mathcal{O}_{\pp \times X}
(\Delta))^{-1} \cdot c(T \pp).$$

\noindent Of course the term $c(T \pp)$ is just $(1+H)^{h^0 -1}$ where
$H$ is the generator of the cohomology of the projective space and
$h^0$ is the rank of $H^0 (X, D)$.  We must therefore understand the
other term in the last 
expression.  Observe

$$\pi_! \mathcal{O}_{\pp \times X} (\Delta) \ = \ R^0 \pi_* - R^1
\pi_* + R^2 \pi_*$$

\noindent in an obvious notation.  But we also have 

$$\mathcal{O}_{\pp \times X}(\Delta) = \pi ^*
\mathcal{O}_{\pp} (1)  \otimes p^* D$$

\noindent and hence $R^i \pi_* \ = \ \mathcal{O}_{\pp} (1) ^{h^i
  (D)}$.  Taking total Chern classes one last time, we deduce that

$$c (\mathrm{Obs}) \ = \ c(\mathcal{O}_{\pp} (1)) ^{-\ind}
\cdot c(T\pp) \ = \ (1+H)^{h^1(D) - h^2(D)}.$$

\noindent Here $\mathrm{ind}$ denotes the index of the
$\cdbar$-operator on the bundle $D$, equivalently the alternating
sum $\sum_i (-1)^i h^i (D)$.  Now we use the assumption that $b_1 (X)
= 0$, and hence $H^1 (X, K_X) \cong H^1 (X; \mathcal{O})^* = 0$
again.  For this implies that 

$$H^0 (X, D) \otimes H^1 (X; K_X - D) \ \rightarrow \ H^1 (X; K_X)$$

\noindent is a trivial map, which in turn implies that whenever $H^0
(X, D) \neq \{0\}$ we have $H^1 (X, K_X-D) = \{0\}$ and hence
$H^1(X,D) = \{0\}$.   Here we use an old result of
Hopf which asserts that whenever $V \otimes V' \rightarrow W$ is a map of
complex vector spaces linear on each factor, then the dimension of the
image exceeds $rk(V) + rk (V') - 1$.
The upshot is that Serre duality implies that the Gromov invariant we
require is exactly

$$Gr_X (\alpha) \ = \ {{-h^2 (D)} \choose {h^0(D) - 1}};$$

\noindent for by the binomial theorem, this is just the Euler class of the
obstruction bundle over the projective space which is the moduli space
of holomorphic curves representing $D$.  Recall the definition of a binomial
coefficient with a negative numerator:

$$ {{-n} \choose {k}} \ = \ (-1)^k {{n+k-1} \choose {k}}.$$

\noindent This shows that

$$Gr_X(\alpha)  =  \pm {{h^0(D)+h^2(D)-2} \choose {h^0(D)-1}} = 
\pm{{h^0(K-D)+h^2(K-D)-2} \choose {h^0(K-D)-1}}  =  \pm Gr_X(\kappa-\alpha).$$

\noindent Here the identity in the middle is given by the Riemann-Roch
theorem and Serre duality for the bundles $\mathcal{O}(D)$ and
$\mathcal{O}(K_X-D)$ on $X$.  If the virtual dimension
$\alpha^2-K_X\cdot\alpha = r$ is in fact
positive, the argument is similar; the relevant equality in this
instance is

$${{-[h^2(D)-r]} \choose {[h^0(D)-1]-r}} \ = \ {{-[h^2(K_X-D)-r]}
  \choose {[h^0(K_X-D)-1]-r}}.$$

\noindent To see this, note that each incidence condition on the
holomorphic curves defines a hyperplane on the projective space of
sections.  The final result again follows from the definition of the
binomial coefficients. 
\end{Pf}

\noindent One can use arguments analogous to those above to perform
explicit computations.  Setting $\alpha = \kappa$ above, we see that
if $b_+ > 1$ then the Gromov
invariant of the canonical class on the surface is $\pm 1$.  This is
the ``holomorphic'' case of the main theorem of \cite{SKDIS}.  Here is a
more substantial assertion.
Suppose $X$ is spin and $c_1 ^2 (X) = 0$; then Rokhlin's theorem shows
that $b_+ (X) = 4n-1$ for some $n \geq 1$.  

\begin{Lem}
For $X$ as above, the Gromov invariant of $D = K_X / 2$ is
  ${{2n-2} \choose {n-1}}$.
\end{Lem}

\noindent To see this, argue as follows. The multiplication map 

$$H^0 (X; D) \otimes H^0 (X; D) \ \rightarrow \ H^0 (X; K_X)$$

\noindent shows that $[b_+ +1]/2 \geq 2h^0 (D)$, by the old result of
Hopf alluded to previously.  But from the Riemann-Roch
theorem, and the duality $H^2 (X;D) \cong H^0 (X;D)^*$ we see that
$2h^0 (X; D) \geq \ind_{\mathcal{O}(D)} (\cdbar) = [b_+ + 1]/2$.
Hence $h^0 (X;D) = [b_+ + 1]/4$.  Now the arguments as above show that 
the obstruction bundle over the projective space $\pp H^0 (X; D)$ satisfies

$$c(\mathrm{Obs}) \ = \ (1+h)^{-\ind} c(T \pp) \ = \ (1+h)^{-[b_+ + 1]/4}.$$

\noindent The Gromov invariant is given by the coefficient of $h^{h^0
  (D) - 1}$ in this expression, which is exactly as claimed by the binomial
  expansion.  (If $X$ is spin but $c_1 ^2 (X) \neq 0$ then the moduli
  space of sections in the class $K_X / 2$ is not zero-dimensional and
  the invariant will vanish as soon as $b_+ > 1$.)  We will return to
  this result in a symplectic setting later, to illustrate a weakness
  of our current theory.


\section{Symplectic surfaces and symmetric products}

In this section we shall review some of the basics of Lefschetz pencils,
and establish the core ideas on which the rest of
the paper is founded.  We include various extensions of results 
from \cite{SKDIS} that will be important later, and have tried to make
the discussion essentially self-contained wherever the statements of
\cite{SKDIS} are inadequate for our applications.  We have also
included an elementary proof of smoothness of the relative Hilbert
scheme for families of curves with only nodal singularities.  Using
this, one can build all the spaces we need - and verify the properties
we need - directly, but at the cost of losing the naturality
which comes for free in the algebraic geometry.  For
foundations on Hilbert schemes and compactified Jacobians, we defer to
the papers \cite{Altman-Kleiman} and \cite{Oda-Seshadri}.


\subsection{Lefschetz pencils}

Let $X$ be an integral symplectic four-manifold: $X$ is
equipped with a symplectic form $\omega$ whose cohomology class
satisfies $[\omega]/2\pi \in H^2(X; \zz) \subset H^2 (X; \rr)$.  Then
there is a line bundle $\mathcal{L}$ with connexion with curvature
$\omega$, and 
Donaldson \cite{Donaldson:submflds} has shown that appropriate sections of high
tensor powers $\mathcal{L}^k$ of 
this line bundle give symplectic surfaces in $X$ Poincar\'e dual to
$k[\omega]/2\pi$.  Extending this in \cite{Donaldson:pencils},
Donaldson shows that 
integral symplectic manifolds admit complex Morse functions or \emph{Lefschetz
pencils}:  we can find a map $f: X \dashrightarrow \sss^2$ defined on
the complement of finitely many points $q_j$ in $X$, with finitely many
critical points $p_i$, all with distinct image under $f$, such that:

\begin{enumerate}
\item $f$ has the local model $(z_1, z_2) \mapsto z_1 z_2$ at each  of
  the $p_i$,
\item $f$ has the local model $(z_1, z_2) \mapsto z_1 / z_2$ at each
  of the $q_j$,
\end{enumerate}

\noindent where all local complex co-ordinates are compatible with
fixed global orientations.  The fibres are again symplectic
submanifolds Poincar\'e dual to the class $k[\omega]/2\pi$, for some
large $k$.  We can assume the symplectic form is positive of
type $(1,1)$ in a fixed almost complex structure at each $p_i, q_j$.
The topology of such a situation is described in \cite{ivanhodge} or
\cite{GompfS} for instance;  after blowing up at each $q_j$ we
have a manifold $X'$ fibred over $\sss^2$ with smooth two-dimensional
surfaces of some fixed genus $g$ as fibres over the points of $\sss^2
\backslash \{f(p_i)\}$ and with critical fibres surfaces with a single
ordinary double point. The exceptional spheres $E_i$ of the blow-ups form
distinct sections of the fibration of $X' = X \sharp r \overline{\cc
  \pp}^2$, where $r = \sharp \{q_j\}$.  We remark that one can always
ensure that the 
singular fibres are irreducible, in the sense that the vanishing
cycles which have collapsed to the node are homologically essential,
or equivalently that removing the node does not separate the singular
fibre \cite{ivanmodulidivisor}.
We will tacitly assume that all Lefschetz fibrations satisfy this
hypothesis henceforth.  The following definition is from \cite{SKDIS}:

\begin{Defn} \label{standardsurface}
Let $f: X' \rightarrow \sss^2$ be a Lefschetz fibration arising from a
Lefschetz pencil $(X,f)$.  A \emph{standard surface} in $X'$ is an
embedded surface $\Sigma \subset (X' \backslash \{ p_i \})$ for which
the restriction $f|_{\Sigma}: \Sigma \rightarrow \sss^2$ is a branched
covering of positive degree with simple branch points.
\end{Defn}

\noindent Here we assume that $f$ is of positive degree on each
component of $\Sigma$ in the case when the surface is disconnected.
As explained in \cite{SKDIS}, a symplectic structure on $X$ induces
a family of symplectic forms $(\omega_{(N)} = p^* (k\omega_X) + N f^*
\omega_{\sss^2})$  on $X'$ which are
symplectomorphic, under the obvious identification away from a small
neighbourhood of
the exceptional sections, to the
original form on $X$ up to scaling by $k(1+N)$.  Moreover, a standard
surface is necessarily symplectic with 
respect to $\omega_{(N)}$ for large enough $N$, and hence a standard
surface disjoint from the $E_i$ gives  rise to a symplectic surface in
$X$.  We will sometimes refer to a finite collection of
standard surfaces with locally positive transverse
intersections as a \emph{positive symplectic divisor}.

\vspace{0.2cm}

\noindent Donaldson's construction of Lefschetz pencils $(X,f)$
involves a subsidiary choice of almost complex structure $J$ on $X$,
compatible with $\omega_X$;  after perturbing $J$ by a $C^k$-small
amount (for any given $k$), we can assume that the fibration $f$ is
pseudoholomorphic and 
all the fibres are almost complex surfaces in $X'$.  Let us also
remark that given any (not necessarily integral) symplectic
four-manifold, arbitrarily small perturbations of the symplectic form
give rise to a rational form and hence to Lefschetz pencils with
fibres dual to a multiple of the form.  Since the $J$-holomorphic
surfaces remain symplectic for small perturbations of $\omega$, this
means that arbitrary symplectic manifolds admit topological  Lefschetz
pencils: that is, if we drop the integrality hypothesis on the form,
we only lose the explicit identification of the homology class of the
fibre.


\subsection{Relative Hilbert schemes}

The choice of $J$ on the total space of the Lefschetz fibration also defines
a smooth map from the base $\phi_f: \sss^2 \rightarrow \mgbar$ to the
Deligne-Mumford moduli space of stable curves, where the extension
over the critical values of $f$ follows precisely from our
requirements on the local normal forms of the singularities.  The map
is defined up to ``admissible isotopy'', that is isotopies which do
not change the geometric intersection number with the divisors of
nodal curves;  these intersections are locally positive.  By choosing
$J$ generically we can assume that the map $\phi_f$ has image disjoint
from the orbifold singular loci of moduli space, and hence lies inside
the fine moduli space of curves without automorphisms.  Hence we have
various universal families over $\phi_f(\sss^2)$, including a
universal curve - which just defines $X'$ - and universal families of
symmetric products and Picard varieties of curves.

\vspace{0.2cm}

\noindent Recall that associated to any smooth Riemann surface
$\Sigma$ we have a
complex torus parametrising line bundles of some fixed degree $r$ and
a smooth complex variety parametrising effective divisors of degree
$r$ on $\Sigma$.  Thinking of these as moduli spaces, for coherent
torsion-free sheaves and zero-dimensional subschemes of fixed length
respectively, machinery from geometric invariant theory
(\cite{moduliofshaves}, \cite{Altman-Kleiman} etc) provides  
relative moduli schemes which fibre over $\mgbar$ - or at least
  over the locus of irreducible curves with at most one node - with fibre
at a smooth point $\tau \in M_g$ just the Picard variety or symmetric
product of $\Sigma_{\tau}$.  These universal fibre bundles give rise
to fibre bundles on the base $\sss^2$ of a Lefschetz pencil, once we
have made a choice of complex structures on the fibres.  More
explicitly, each of these fibre bundles can also be defined by
local charts, called \emph{restricted charts} in \cite{SKDIS}.  A
restricted chart is a diffeomorphism $\chi: D_1 \times D_2 \rightarrow X'$
which is a smooth family, indexed by $\tau \in D_1 \subset \sss^2$, of
holomorphic diffeomorphisms $\chi_{\tau}: \{ \tau \} \times D_2 \rightarrow
U_{\tau} \subset f^{-1} (\tau)$ onto open subsets of the
fibres of $X'$.  Here each $D_i$ denotes the unit disc in the complex
plane, with its standard integrable structure.  The existences of
atlases of restricted charts is implied by the Riemann mapping
theorem with smooth dependence on parameters.  At any rate, with this
background we can make the following: 

\begin{Defn} \label{Hilbetc}
Let $(X,f)$ be a Lefschetz pencil inducing $\phi_f: \sss^2 \rightarrow
\mgbar$.  
\begin{itemize}
\item Denote by $F:X_r(f)\rightarrow \sss^2$ the pullback by $\phi_f$ 
of the universal relative Hilbert scheme for zero-dimensional length
$r$ subschemes of fibres of the universal curve $\pi: \mathcal{C}_g
\rightarrow \mgbar$.  
\item Denote by $G:P_r (f)\rightarrow \sss^2$ the pullback by $\phi_f$
of the universal relative Picard scheme for degree $r$ torsion-free
sheaves on the fibres of $\pi: \mathcal{C}_g \rightarrow \mgbar$.
\end{itemize}
\noindent The natural smooth map $u: X_r(f) \rightarrow
P_r(f)$ will be referred to as the \emph{Abel-Jacobi map}.
\end{Defn}

\noindent The existence of the map $u$ is proven in
\cite{Altman-Kleiman} or can be deduced (in the smooth category) from
our constructions.  Here is a set-theoretic description of the
singular fibres of $G$ and $F$.

\begin{enumerate}
\item The degree $r$ torsion-free sheaves on an irreducible nodal
  curve $C_0$ are of 
  two forms: locally free, or push-forwards of locally free sheaves
  from the normalisation $\pi: \tilde{C}_0 \rightarrow C_0$.  A
  locally free sheaf is completely
  determined by a pair $(L, \lambda)$ where $L \rightarrow
  \tilde{C}_0$ is a degree $r$ line bundle on the normalisation and
  $\lambda \in \mathrm{Iso}(L_{\alpha}, L_{\beta}) \cong \cc^*$ is a
  gluing paramater which identifies the fibres $L_{\alpha}$ and
  $L_{\beta}$ of $L$ over the preimages of the node of $C_0$.  This
  gives a $\cc^*$-bundle over $\Pic_r(\tilde{C}_0)$.  The
  non-locally free sheaves are of the form $\pi_* L'$ where $L'
  \rightarrow \tilde{C}_0$ is locally free of degree $r-1$.  These
  arise by compactifying $\cc^*$ to $\pp^1$ and identifying the two
  degenerate gluings - the $0$ and $\infty$ sections of the resulting
  $\pp^1$-bundle  - over a
  translation by the action of $\mathcal{O}(\alpha - \beta)$ in
  $\Pic_r(\tilde{C}_0)$. 

\item The fibre of $F$ can be completely described by giving the fibre
  of $G$, as above, and the fibres of the map $u$.  The latter are
  projective spaces.  At a point $(L, \lambda)$ the fibre of $u$ is
  the subspace of the linear system $\pp H^0 (L)$ comprising those
  sections $s \in H^0 (L)$ for which $s(\alpha) = \lambda s(\beta)$.
  At a point $\pi_* L'$ the fibre of $u$ is just the entire linear
  system $\pp H^0 (L')$.
\end{enumerate}

\noindent The sets $X_r(f)$ and $P_r(f)$ obtained above carry obvious
topologies:  given a sequence $(D_n = p_n + D)_{n \in \nn}$ of distinct
$r$-tuples of points in the singular fibre of $f$, with $D$ a fixed
$(r-1)$-tuple and $p_n \rightarrow \mathrm{Node}$, then the points of
the Hilbert scheme converge to the obvious point of
$\Sym^{r-1}(C_0)$ which is determined by the divisor $D$ and the
associated line bundle $L' = \mathcal{O}(D)$ of rank $r-1$ on
$\tilde{C_0}$.  The general behaviour is analogous.  We can put smooth
structures on the spaces using explicit local charts.
The following result - which may be known to
algebraic geometers but does not appear in the literature - is central
for this paper.   

\begin{Thm}
For any $(X,f)$ and each $r$, the total spaces of $X_r(f)$ and
$P_r(f)$ are smooth compact symplectic manifolds.
\end{Thm}

\begin{Pf}
Up to diffeomorphism, the total space of the relative Picard variety
for a fibration with a section is independent of $r$; 
hence the proof given for $r=2g-2$ in \cite{SKDIS} is sufficient.
We recall the main point: any torsion-free sheaf of degree $r$ on an
irreducible 
nodal curve is either locally free or is the push-forward of a locally
free sheaf of degree $r-1$ on the normalisation.  This follows from
the existence of a short exact sequence 

$$0 \rightarrow \mathcal{O}_C \rightarrow \pi_*
\mathcal{O}_{\tilde{C}} \rightarrow \cc_{(p)} \rightarrow 0$$

\noindent where $\pi: \tilde{C} \rightarrow C$ is the normalisation
and the skyscraper sheaf $\cc_{(p)}$ is supported at the node.  The
locally free sheaves on $C$ come from degree $r$ locally free sheaves
on $\tilde{C}$ with a gluing parameter $\lambda \in \cc^*$ to identify
the fibres of the line bundle at the two preimages of the node.  The
compactification by adding torsion free sheaves arises from
compactifying $\cc^*$ (as a $\cc^* \times \cc^*$-space) to $\cc \pp^1$ and
gluing together the $0$- and $\infty$-sections over a translation in
the base, cf. the Appendix to \cite{SKDIS}.  The resulting variety has
normal crossings and the total space of the relative Picard scheme is
easily checked to be smooth, modelled transverse to the singularities
of the central fibre on a family of semistable elliptic curves. 

\vspace{0.2cm}

\noindent  For the relative Hilbert scheme, we can argue as follows.
Clearly we 
have smoothness away from the singular fibres, and more generally when
the subscheme is supported away from the nodes of fibres of $f$.
Moreover, given a
zero-dimensional subscheme of a nodal fibre $\Sigma_0$ of $f$
supported at a collection of points $x_1, \ldots, x_s$ we can take
product charts around each of the $x_i$ and reduce to the situation at
which all of the points lie at the node of $\Sigma_0$.  Hence it will
be sufficient to prove smoothness for the local model: that is, for
the relative Hilbert scheme $\cc_r(f)$ of the map $f: \cc^2
\rightarrow \cc$ defined by $(z,w) \mapsto zw$.

\vspace{0.2cm}

\noindent According to Nakajima \cite{Nakajima}, the Hilbert scheme
$\Hilb^{[r]}(\cc^2)$ is globally smooth and may be described
explicitly as follows.  Let $\tilde{\mathcal{H}}$ denote the space

$$\big \{ (B_1, B_2, v) \in M_r(\cc)^2 \times \cc^r \ | \ [B_1, B_2] =
0, \ (*) \big \}$$

\noindent where the stability condition $(*)$ asserts that for any $S
\subsetneq \cc^r$ invariant under both $B_i$ we have $v \not \in S$.
There is a $GL_r(\cc)$ action on $\tilde{\mathcal{H}}$ defined by

$$g \cdot (B_1, B_2, v) \ \mapsto \ (gB_1g^{-1}, gB_2g^{-1}, gv)$$

\noindent and the stability condition implies that this is free.
Moreover, the cokernel of the map 

$$\tilde{\mathcal{H}} \ \rightarrow \ M_r(\cc); \qquad (B_1, B_2, v)
\mapsto [B_1, B_2]$$

\noindent has constant rank $r$ (the stability condition shows that
the map $\xi \mapsto \xi(v)$ is an isomorphism from the cokernel to
$\cc^r$). It follows that $\tilde{\mathcal{H}}$ is
smooth, and the freeness of the action gives a smooth structure on the
quotient $\mathcal{H} = \tilde{\mathcal{H}}/GL_r(\cc)$; but this
quotient is exactly $\Hilb^{[r]}(\cc^2)$.  To see this, note that an
ideal $\curly{I} \subset \cc[z,w]$ of length $r$ defines a quotient
vector space  $V = \cc[z,w] / \curly{I}$ of dimension $r$; now define
endomorphisms $B_1, B_2$ via 
multiplication by $z,w$ respectively, and $v$ by the image of $1$.
Conversely, given $(B_1, B_2, v)$ we set $\curly{I} = \ker (\phi)$ for
$\phi: \cc[z,w] \rightarrow \cc^r$ defined by $f \mapsto f(B_1,
B_2)v$; this is the inverse map.  There is an obvious family of
zero-dimensional subschemes over $\mathcal{H}$; given any other such
family $\pi:Z \rightarrow W$, there is a locally free sheaf of rank
$r$ over $W$, given by $\pi_* \mathcal{O}_Z$.  Taking a cover of $W$
and trivialising locally, we can define multiplication maps $B_i$ by the
co-ordinate functions $z_i$ and hence show this family is indeed a
pullback by a map $W \rightarrow \mathcal{H}$, which gives the
required universal property.

\vspace{0.2cm}

\noindent This gives a simple description of the relative Hilbert scheme
$\mathcal{H}(f)$ as the subscheme of those ideals containing
$zw-\lambda$ for some $\lambda \in \cc$.  Explicitly, setting $I_r$
for the identity element of $End(\cc^r)$, we have:

$$\tilde{\mathcal{H}}(f) \ = \ \big \{ (B_1, B_2, \lambda, v) \in
M_r(\cc)^2 \times \cc \times \cc^r \ | \ B_1 B_2 = \lambda I_r, \,
B_2 B_1 = \lambda I_r, \, (*) \big \};$$

\noindent this still carries a free $GL_r(\cc)$ action (trivial on the
$\lambda$-component) and the quotient is the relative Hilbert scheme
$\mathcal{H}(f)$.  Note that if $\lambda \neq 0$ then each of the
$B_i$ is invertible, the first two equations are equivalent, and the
two matrices are simultaneously diagonalisable; then the fibre over
$\lambda \in \cc^*$ is just the obvious copy of $\cc^r =
\Sym^r(\cc)$.  Smoothness of the total space now follows from the fact
that the map

$$M_r(\cc)^2 \oplus \cc \ \rightarrow \ M_r(\cc)^2; \qquad (C_1, C_2,
\mu) \mapsto (C_1 B_2 + B_1 C_2 - \mu I_r, B_2 C_1 + C_2 B_1 - \mu
I_r)$$

\noindent (which is the differential of the defining equations) has
constant dimensional kernel $r^2+1$, independent of $(B_1, B_2,
\lambda, v) \in \tilde{\mathcal{H}}(f)$.  We can see this with a
routine if tedious computation.  First, if $\lambda \neq 0$, then the
equation $B_1 B_2 = \lambda I_r$ shows that each $B_i$ is invertible;
the two defining equations are then equivalent, and at any point of
the kernel, $C_1$ is uniquely determined by $C_2$ and $\mu$ which can
be prescribed arbitrarily:

$$\mathrm{Kernel} \ = \ \{ (C_1, C_2, \mu) \ | \ C_1 = \frac{\mu B_1 -
  B_1 C_2 B_1}{\lambda} \}.$$

\noindent Now suppose $\lambda = 0$ so we are in the singular
fibre. Suppose in addition that neither $B_1$ nor $B_2$ vanish
identically.  The stability condition $(*)$ implies that we can choose
a basis for $\cc^r$ of the form

$$\cc^r \ = \ \langle v, B_1 v, B_1^2v, \ldots B_1^{n}v, B_2v,
B_2^2v, \ldots, B_2^m v \rangle; \qquad r=n+m+1.$$

\noindent For certainly the space $Span \langle B_1^iv, B_2^jv
\rangle_{i,j\geq 0}$ is $B_i$-invariant and contains $v$, so must be
full, and $v \neq 0$ or $\{ 0 \} \subset \cc^r$ violates $(*)$.
Suppose first that neither $B_1$ nor $B_2$ identically vanishes.  Then
we can write the $B_i$ in matrix form, with respect to the above
basis, as follows: write $O_{a,b}$ for the $(a \times b)$-matrix with
all entries $0$.

$$
B_1 \ = \ \left( 
\begin{tabular}{c|ccccc|c}
$0$   &    $0$ & $0$ & $\cdots$ & $0$ & $h_0$      & $0_{1,m}$  \\ \hline
$1$   &    $0$ & $0$ & $\cdots$ & $0$ & $h_1$      &  $\vdots$        \\
$0$   &    $1$ & $0$ & $\cdots$ & $0$ & $h_2$      &                \\
$0$   &    $0$ & $1$ & $\cdots$ & $0$ & $h_3$      &   $0_{n,m}$ \\
$\vdots$&      &   & $\vdots$ &   &          &                 \\
$0$     &  $0$ & $0$ & $\cdots$ & $1$ & $h_n$      &    $\vdots$  \\ \hline
$0_{m,1}$&   $\ldots$ & & $0_{m,n}$ & & $\ldots$   &   $0_{m,m}$        
\end{tabular}
\right); 
$$

\vspace{0.2cm}

$$ 
B_2 \ = \ \left(
\begin{tabular}{c|c|ccccc}
$0$    &    $0_{1,n}$   &   $0$  &  $0$  & $\cdots$  &  $0$  & $H_0$ \\ \hline
$0_{n,1}$& $0_{n,n}$   &      &     &  $0_{n,m}$ &     &    \\ \hline
$1$   &    $\vdots$    &  $0$  &  $0$  & $\cdots$  &  $0$  & $H_1$ \\
$0$     &              &   $1$  &  $0$  & $\cdots$  &  $0$  & $H_2$ \\
$0$     &   $0_{m,n}$   &   $0$  &  $1$  & $\cdots$  &  $0$  & $H_3$ \\
$\vdots$&              &      &     & $\vdots$  &     &     \\
$0$     &    $\vdots$    &   $0$  &  $0$  & $\cdots$  &  $1$  & $H_m$ 
\end{tabular}
\right) 
$$    

\noindent Here $B_1^{n+1}v = \sum h_i B_1^i v$, similarly $B_2^{m+1}v
= \sum H_i B_2^i v$, using the fact that the subspaces $\langle v,
B_jv, B_j^2v, \ldots \rangle$ are invariant under $B_j$, $j = 1,2$.
It is now easy to check that $h_0 = 0 = H_0$, from the condition that
$B_1B_2 = B_2B_1 = 0$ and neither $B_j \equiv 0$.  Let us then write,
in obvious notation, an element $(C_1, C_2, \mu)$ in the kernel of the
differential as follows, in block matrix form for the $C_i$:

$$
C_1 \ = \ \left( 
\begin{tabular}{c|c|c}
$\beta$         & $\tau$       & $\chi$      \\ \hline
$\tilde{\tau}$  & $\mu_{i,j}$  & $\nu_{i,j}$    \\ \hline
$\tilde{\chi}$& $\phi_{i,j}$& $\psi_{i,j}$   
\end{tabular}
\right); \qquad
C_2 \ = \ \left(
\begin{tabular}{c|c|c}
$\alpha$        & $v$           & $u$           \\ \hline
$\tilde{v}$     & $m_{i,j}$     & $n_{i,j}$     \\ \hline
$\tilde{u}$     & $p_{i,j}$     & $q_{i,j}$     
\end{tabular}
\right)
$$

\noindent Here $\tau, v$ are $(1\times n)$-vectors, $\chi, u$ are
$(1 \times m)$-vectors (and similarly for the
$\tilde{\cdot}$-entries); whilst $m_{i,j}, \mu_{i,j}$ are in
$M_n(\cc)$, $q_{i,j}, \psi_{i,j} \in M_m(\cc)$ and the other entries
are $(n\times m)$ or $(m \times n)$ blocks in the obvious way.  In
this schematic, the linearisation of the defining equations for our
relative Hilbert scheme become:

$$
\mu I_r \ = \ B_2 C_1 + C_2 B_1 \ = \ \left(
\begin{tabular}{c|c|c}
$v_{1,1}$ & $\Phi v$ & $0$ \\ \hline
$m_{i,1}$ & $(m_{i,j})\Phi$ & $0$ \\ \hline
$\beta + \Psi \tilde{\chi} + p_{j,1}$ & $\tau + \Psi (\phi_{i,j}) +
(p_{i,j})\Phi$  & $\chi + \Psi (\psi_{i,j})$
\end{tabular} \right)
$$

\noindent and simultaneously the equation:

$$
\mu I_r \ = \ B_1 C_2 + C_1 B_2 \ = \ \left(
\begin{tabular}{c|c|c}
$\chi_1$ & $0$ & $\chi'$ \\ \hline
$\alpha + v_{1,1} + (\tilde{v})'$ & $v^t + \Phi (m_{i,j})$ & $u^t +
\Phi(n_{i,j}) + (\nu_{i,j}) \Psi$ \\ \hline
$\psi_{k,1}$ & $0$ & $(\psi_{i,j}) \Psi$
\end{tabular}
\right) 
$$

\noindent Here we have written $\Phi$ and $\Psi$ for the non-trivial
$(n\times n)$, resp. $(m \times m)$, blocks in the matrices $B_1$ and
$B_2$.  In addition, $\chi'$ is the vector $(\chi_2, \ldots,
\chi_{m-1}, \sum H_{i+1} \chi_i)$ and $(\tilde{v})'$ is
similarly a linear combination of the entries of $\tilde{v}$ and the
$h_i$ (the precise formula is not really important, and easily worked
out by the curious).
First of all, we claim that necessarily $\mu = 0$; i.e. the
matrices above must vanish.  For from the first equation, we see from
the left hand column that $m_{i,1}$ vanishes for each $i$, and hence
the matrix $(m_{i,j})$ has trivial determinant; but then we cannot
have $(m_{i,j}) \Phi = \mu I_n$ unless $\mu = 0$.  Hence the equations
for our differential amount to the vanishing of the equations above.
Clearly we can freely choose the parameters $q_{i,j}, \mu_{i,j}, \tilde{u}$ and
$\tilde{\tau}$, which do not even appear on the RHS.  A moment's
inspection shows that one can also prescribe $\alpha, \beta, u, \tau$
freely, and then $p_{i,j}$ and $n_{i,j}$; then all the other data is
determined.  This is clear if $\Phi$ and $\Psi$ are invertible, just
by manipulating the above.  In this case, an element of the kernel is
completely determined by the choice of

$$\alpha, \beta, u, \tau, \tilde{u}, \tilde{\tau}, q_{i,j}, \mu_{i,j},
p_{i,j}, n_{i,j}$$

\noindent of respective dimensions:

$$1,1,m,n,m,n,m^2,n^2,mn,mn.$$

\noindent Hence the dimension of the kernel is $(m+n+1)^2+1 = r^2+1$,
as required.  We claim this is still the case even if $\Phi$ and
$\Psi$ are not necessarily invertible.  For instance, we already know
that the matrix $m_{i,j}$ is of the form

$$m_{i,j} \ = \ \big( 0_{n,1} \ | \ m_{i,j}' \big)$$

\noindent for some $(n \times (n-1))$-matrix $m_{i,j}'$. Then if in
fact $\Phi$ has trivial determinant, we can see that 

$$\Phi \ = \ \left( \begin{tabular}{cc} $0_{n-1,1}$&$0$ \\ \hline
    $I_{n-1}$ & $(h_j)$ 
  \end{tabular} \right)$$

\noindent from which it easily follows that $(m_{i,j})\Phi = 0
\Rightarrow (m_{i,j}) \equiv 0$.  Then the rest of the data can be
determined successively.

\vspace{0.2cm}

\noindent This leaves only the case where (without loss of generality)
$B_1 \equiv 0$.  In this case, it is no longer true that $\mu$
necessarily vanishes; but it is easy to check the kernel is now
$(r^2+1)$-dimensional.  For instance, if $\mu \neq 0$ then $\mu$
determines $C_2 = \mu B_1^{-1}$ and $C_1$ can be chosen freely.  It
follows that the space $\tilde{\mathcal{H}}(f)$ is determined
everywhere by a map to a vector space of constant rank, hence is
smooth, and then the quotient by the free $GL_r(\cc)$ action is smooth
by Luna's slice theorem.  Hence we have smooth compact manifolds.  To
put symplectic structures on the spaces $X_r(f)$ follows from a
theorem of Gompf \cite{GompfGokova}, adapting an old argument of
Thurston, and is discussed below. 
\end{Pf}
\begin{Rmk}
One can use the above to \emph{define} the relative Hilbert scheme, as
a smooth manifold via local charts, for readers wary of sheaf
quotients.  Note also that the above shows that we have a smooth
compactification of the symmetric product fibration even when there
are reducible fibres present.
\end{Rmk}

\noindent In fact the spaces $X_r(f)$ all
have \emph{normal 
  crossing singularities} (they are locally of the form $\{ zw=0 \}
\times \{ \mathrm{Smooth} \}$).  For large $r$ this follows from the
  result for the relative Picard fibration and the existence of the
  Abel-Jacobi map, whose total space is a projective bundle over the
  base once $r>2g-2$.  Since the condition is local, the result
  for large $r$ can be used to deduce it in general.  The $r$-th
  symmetric product of 
$\{ zw=0 \}$ comprises the spaces $\cc^{r-s} \times \cc^s$ glued
together along various affine hyperplanes, and the Hilbert scheme is
obtained by successively blowing up the strata; the projective
co-ordinates in the normalisation are described in \cite{Nakajima} and
  the Appendix to \cite{SKDIS}.  A ``conceptual'' proof of smoothness
  from this point of 
  view is also given there, though the computations above apply more
  easily to the case of small $r$.

\vspace{0.2cm}

\noindent Since we have integrable complex structures near the
singular fibres, we can extend (for instance via a connexion) to
obtain almost complex structures on the total spaces of the
$X_r(f)$. For us, the holomorphic structure on the fibres is always standard.
The existence of (a canonical deformation equivalence
class of) symplectic structures on $X_r(f)$ runs as follows.  The
obvious integrable K\"ahler forms on the fibres patch to give a
global two-form $\Omega_0$; then if $\omega_{st}$ is the area form on
the sphere, with total area $1$, the forms $\Omega_t = \Omega_0 + t
F^* \omega_{st}$ are symplectic for all sufficiently large $t$. 

\vspace{0.2cm}

\noindent Write $\mathcal{J} = \mathcal{J}_r$ for the class of almost
complex structures $J$ on $X_r(f)$ which agree 
with the standard integrable structures on all the fibres and which
agree with an 
integrable structure induced from an almost complex structure on $X'$
in some tubular neighbourhood of each critical fibre.
Usually we will deal with smooth almost complex structures, though in
the next section it shall be important to allow ones that are only
H\"older continuous.  In any case, all our almost complex structures
shall be of this form.

\begin{Lem} \label{connectedparameters}
The space $\mathcal{J}$ is non-empty and contractible. A given $J \in
\mathcal{J}$ is tamed by the symplectic forms $\Omega_t$ for all
sufficiently large $t > t(J) \in \rr_+$.
\end{Lem}

\noindent This is just as in \cite{SKDIS}; if we fix the
neighbourhoods of the critical fibres where the structure is induced
by an integrable structure on $X$ then the remaining choice is of a
section of a fibre bundle with contractible (affine) fibres.  Note
that the complex structures $J$ will not be compatible with
$\Omega_t$;  the class of compatible structures is not large enough to
achieve regularity for spaces of holomorphic sections.

\vspace{0.2cm}

\noindent Let us recall various facts
concerning the geometry of the Abel-Jacobi map.  As we observed at the
start of the section, the fibres of $u$ are projective spaces
(linear systems).  We can check their generic dimension at the
singular fibre.  A locally free sheaf
$L$ on $C_0$ is given by a line bundle $\tilde{L}$ on
$\tilde{C}_0$ with a
gluing map $\lambda: \cc \rightarrow \cc$ of the fibres over the two
preimages $p,q$ of the node.  Since $\tilde{C}_0$ is of genus $g-1$ when
$C_0$ has (arithmetic) genus $g$, if $L$ has degree $d$ then by the
Riemann-Roch theorem we find that $\ind_{\tilde{L}} (\cdbar) =
d-(g-1)+1$.  Generically this is the dimension of $H^0 (\tilde{L})$
and then 

$$H^0 (L) = \{ s \in H^0 (\tilde{L}) \ | \ s(p) = \lambda s(q) \}.$$

\noindent Provided $\lambda \in \cc^*$ and not every section vanishes
at $p,q$ this is a hyperplane; hence $h^0 (L) = d-g+1$.  Thus the
generic dimension of the fibres of the Abel-Jacobi map is the same
over the singular fibre of the Picard fibration as over the smooth
fibres.  There are two different behaviours we should emphasise:

\begin{itemize}
\item Along the normal
crossing divisor of the singular fibre of $P_r(f)$ we have $\lambda
\in \{0, \infty\}$ (these two cases are glued together,
cf. \cite{SKDIS}).  Here all sections
represent Weil divisors which are not Cartier and do not arise from
locally free sheaves; for instance, any subscheme supported at the
node to order exactly one is of this form.  
\item We also have an embedding (of varieties
of locally free sheaves)
$\Pic_{d-2}(C_0) \rightarrow \Pic_d (C_0)$ induced by the embedding
$\Pic_{d-2}(\tilde{C}_0) \rightarrow \Pic_d (\tilde{C}_0)$ which takes
a line bundle $\tilde{L} \rightarrow \tilde{L} \otimes
\mathcal{O}(p+q)$.  At these points,
the space of sections of $H^0 (L)$ is the entirety of the space of
sections of $H^0 (\tilde{L} \otimes \mathcal{O}(-p-q))$ - the parameter
$\lambda$ plays no role now - which although
not a hyperplane again has the right dimension, since the bundle
upstairs has a different degree.  
\end{itemize}

\noindent In particular, surfaces in the four-manifold $X'$ which pass
through a node of a singular fibre can arise from smooth sections of the
symmetric product fibration (necessarily disjoint from the locus of
critical values for the projection $F$).  If this was not the case, we
could not hope 
to obtain a compactness theorem for spaces of sections.

\vspace{0.2cm}

\noindent In \cite{SKDIS}
we didn't identify the monodromy of these associated families of
Jacobians or symmetric products, so let us do that here.    
Remark for the readers' convenience that this shall not be used later
on and is included just for the sake of intuition.  The author learned
this general material from unpublished notes of Paul Seidel
\cite{Seidel:thesisdraft}.

\begin{Defn} \label{gendehntwist}
Let $(X^{2n},\omega_X)$ be a symplectic manifold.  Be given a symplectic
manifold $(Y^{2n-2},\omega_Y)$ 
and a principal circle bundle $\sss^1 \rightarrow W
\stackrel{\pi}{\longrightarrow} Y$, together with an embedding $\iota:W
\hookrightarrow X$ which satisfies $\iota^* \omega_X = \pi^*
\omega_Y$.  Then a \emph{generalised Dehn twist} along $W$ is a
symplectomorphism $\phi: X \rightarrow X$ with the following two
properties:

\begin{itemize}
\item $\phi$ is the identity outside a tubular neighbourhood of $W$;
\item on each circle fibre $\phi$ acts as the antipodal map.
\end{itemize}

\noindent (If $n=1$ and $Y$ is a point then this is a Dehn twist along
a curve in a real surface in the classical sense.) 
\end{Defn}

\noindent It is not hard to show that such a symplectomorphism $\phi$
always exists in this situation (and indeed for more general
coisotropic embeddings of sphere bundles with orthogonal structure
group). An important result, due to Seidel, is that the data
distinguishes a unique Hamiltonian isotopy class of
symplectomorphism containing generalised Dehn twists.  With this
terminology established, we then have: 

\begin{Prop}
Let $X$ be a family of genus $g$ curves over the disc $D$ with a
single irreducible nodal fibre $\Sigma_0$ over $0$.  Write
$\tilde{\Sigma}_0$ for the normalisation of this nodal fibre.
\begin{itemize}
\item The monodromy of the Picard fibration $P_r(f) \rightarrow D$
  around $\partial D$ is a generalised Dehn twist along a circle
  bundle over an embedded copy of $\Pic_{r-1} (\tilde{\Sigma}_0)$.
\item The monodromy of the relative Hilbert scheme $X_r(f) \rightarrow
  X$ around $\partial D$ is a generalised Dehn twist along a circle
  bundle over an embedded copy of $\Sym^{r-1} (\tilde{\Sigma}_0)$.
\end{itemize}

\noindent In each case, the base of the circle bundle is
isomorphic to the singular locus of the fibre over $0$, and the
singularity arises as the circle fibres shrink to zero size as we
approach the origin of the disc.
\end{Prop}

\begin{Pf}[Sketch]  The point is that for any
  holomorphic fibration with a normal crossing fibre over $0 \in D$,
  the monodromy 
  is a generalised Dehn twist of this form.  This is because the
  normal crossing data defines a unique local model in a neighbourhood
  of the singular locus:  the circle bundle is just the bundle of
  vanishing cycles defined with respect to some symplectic connexion
  arising from a local K\"ahler form.  Seidel's notes give a detailed
  construction of the symplectomorphism from this data.  It follows
  that it is enough to identify the singular locus of the singular fibres. 
(Of course this discussion holds for $X$ itself.  If we fix a smooth fibre
over a point $t \in D^*$, an embedded ray $\gamma$ from $t$ to $0$
  in $D$ and a symplectic
form on the total space of $X$, then the fibre $X_t$ contains a
distinguished real circle, which is the vanishing cycle associated to
the critical point of $f$ by the path $\gamma$. The monodromy is a
  Dehn twist about this circle.)

\vspace{0.2cm}

\noindent $\bullet$ The Picard fibration $W$, up to diffeomorphism
(for instance on choosing a local section), is just the total space of
the quotient $R^1f_*(\mathcal{O}) / R^1f_*(\zz)$.  
The singular locus of $W_0$ is just a copy of the Jacobian of the
normalisation of the singular fibre $X_0$:  

$$Sing(W_0) \ \cong \ \Jac(\tilde{X}_0) \ \cong \ \rot^{2g-2}.$$

\noindent This follows from the explicit construction of the
generalised Jacobian given in the Appendix to \cite{SKDIS}.
In fact we can see clearly why
the monodromy of the fibration is a generalised Dehn twist.  From our
smooth description, the total space is
diffeomorphic to $\Jac(\tilde{X}_0) \times E$, with $E$ a fibration
of elliptic curves with a unique semistable nodal fibre.  This is
because the integral homology lattices give a flat connexion on a real
codimension two subbundle of the homology bundles.  The monodromy
is just $\mathrm{id} \times \tau_{\gamma}$ where $\tau$ is the
diffeomorphism induced by (the homological action of) the Dehn
twist about $\gamma$ on the subspace of $H_2(X_t)$ generated by $\gamma$
and a transverse longitudinal curve. 
Invariantly, the singular locus is the degree $r-1$ Picard variety of
the normalisation since the natural map $\mathcal{O}_{X_0} \rightarrow
\pi_* \mathcal{O}_{\tilde{X}_0}$ has cokernel of length one.

\vspace{0.2cm}

\noindent $\bullet$ Consider now the
associated family $Z = \Hilb^{[r]}(f)$.  When $r$ is large (say
$r > 2g-2$), we can use 
the above and the Abel-Jacobi map.  For then the Hilbert scheme is a
projective bundle over the Picard variety, and the singular locus of
the zero fibre is a projective bundle over $\Pic_{r-1}
(\tilde{\Sigma}_0)$.  The total space of
this bundle is isomorphic to a copy of $\Sym^{r-1}
(\tilde{\Sigma}_0)$.  Naively
this description breaks down for small $r$, although the
result is still true;  it is a consequence, for $r=2$, of the examples
described in the Appendix to \cite{SKDIS}, and more generally of the
fact that the singular locus represents precisely Weil non-Cartier
divisors.  Now the fact that every torsion-free sheaf is either
locally free or the push-forward from the normalisation of something
locally free implies the Proposition; sections of sheaves $\pi_*
\mathcal{L}$ correspond to the obvious linear systems on the
normalisation. 
\end{Pf}


\section{Counting standard surfaces}

In this section, we present the holomorphic curve theory for our
associated fibrations.  We shall compute the virtual dimensions of
moduli spaces, prove that in fact they are \emph{a priori} compact for
regular almost complex structures, and
define the appropriate Gromov invariant counting sections.


\subsection{Sections, cycles and index theory}

From a Lefschetz pencil $(X,f)$ we build the fibration $X'$ and the
associated fibrations $F = F_r: X_r(f) \rightarrow \sss^2$.
Suppose $s: \sss^2 \rightarrow X_r(f)$ is a smooth section.
At every point $t \in \sss^2$ we have a collection of $r$ points on
the fibre $\Sigma_t$ of $f: X' \rightarrow \sss^2$.  As $t$ varies
these trace out some cycle $C_s \subset X'$ and the association $s
\mapsto [C_s]$ defines a map $\phi: \Gamma(F) \rightarrow H_2 (X';
\zz)$ from homotopy classes of sections of $F$ to homology classes of
cycles in $X'$.  

\begin{Lem} \label{Cyclesinject}
The map $\phi: \amalg_r \Gamma(F_r) \rightarrow H_2 (X';\zz)$ is
injective. 
\end{Lem}

\begin{Pf}
Clearly the image homology class determines $r$, as the algebraic
intersection number with the fibre of $f$.  So we can restrict to a
single $F=F_r$ henceforth.
We will construct a ``partial inverse''.  Fix $A \in H_2 (X'; \zz)$;
this defines a complex line bundle $L_A \rightarrow X'$.  A fixed
symplectic form on $X'$ defines a connexion (field of horizontal
subspaces) on the smooth locus of $X' \rightarrow \sss^2$; also pick a
connexion $\nabla$ on the total space of the principal circle bundle
$\mathcal{P}_A \rightarrow X'$.  Now consider the circle bundle
$\mathcal{P}_t$ given by restricting $\mathcal{P}_A$ to a smooth fibre
$C_t = f^{-1}(t)$ of the Lefschetz fibration.  At each point $u \in
\mathcal{P}_t$ we have a natural decomposition of the tangent space of
$(\mathcal{P}_A)_u$ given by 

$$T_u (\mathcal{P}_A) \ = \ T_u (\sss^1) \oplus T_u (C_t) \oplus f^*
T_{f(u)} (\sss^2)$$

\noindent using first the connexion $\nabla$ and then the symplectic
form.  In particular, the tangent space $T_u (\mathcal{P}_t)$ is
naturally split and this splitting is obviously invariant under the
action of $\sss^1$; hence we induce a natural connexion on each of the
complex line bundles $L_t = L_A |_{C_t}$.  But every connexion induces an
integrable holomorphic structure over a one-dimensional complex
manifold (i.e. for each $L_t$ we have $\cdbar^2 = 0$).  Thus the
(contractible) choice of connexion $\nabla$ gives rise to a
distinguished section of the bundle $P_r(f)$ which is just the class
of the holomorphic line bundle $L_t \rightarrow C_t$ as $t$ varies
over $\sss^2$.  Riemann's removable of singularities theorem, with
holomorphic data near the singular fibres, takes care of the extension
even where the connexion becomes singular.  In other words, $H_2 (X';
\zz)$ is in one-to-one
correspondence with the set of different homotopy classes of sections
of bundles $P_r (f)$ as $r$ varies.  (One can also see this by thinking
of connexions on $L_t$ as an infinite dimensional affine space;
any fibre bundle of such spaces admits a section.)

\vspace{0.2cm}

\noindent Now a section of $P_r(f)$ may not
lift to any section of $X_r(f)$, but if it does lift then the
different homotopy classes of lifts differ by at most an action of
$\zz$ which acts by adding a generator $h$ of $\pi_2 (\pp^N)$ to a given
section: here $\pp^N$ is a linear system fibre of the Abel-Jacobi
map. (These projective spaces may be empty, or have varying dimension
as we move the section by homotopy, but they contribute at most a
single integral class to $\pi_2$ \cite{ACGH}.) Now it is easy to check
directly that if a section $s$ of
$X_r(f)$ defines a homology class $A_s \in H_2 (X';\zz)$ then $s + h$
defines a homology class $A_s + [\mathrm{Fibre}]$, and these of course differ.
It follows that at most a single homotopy class of section can yield
any given $A$ under the above map $\phi$.
\end{Pf}

\noindent We will also need to introduce the \emph{twisting} map $\iota: H_2
(X; \zz) \rightarrow 
H_2 (X';\zz)$ which takes a cycle $C \subset X$ to the cycle $C
\cup E$ where $E$ is the union of the exceptional sections of $f$,
each taken with multiplicity one.  Clearly this is an embedding of
$H_2 (X; \zz)$ onto a direct summand of $H_2 (X'; \zz)$. We will
sometimes identify 
homology classes on $X$ with elements in the homology of $X'$ under
the natural map $i$ on $H_2$ induced by the blow-down map; that is, use
Poincar\'e duality and the pullback on $H^2$, or just choose a cycle
not passing through the points we blow-up and take its preimage.  Note
this ``obvious'' embedding $i:H_2 (X;\zz) \mapsto H_2(X';\zz)$ is
different from the twisting map: $\iota(\alpha) = i(\alpha) + \sum E_i$.
With this convention, for any $\alpha \in H_2
(X; \zz)$ there is an equality

\begin{Eqn} \label{samedim}
\alpha^2 - \alpha \cdot K_X \ = \ \iota(\alpha)^2 - \iota(\alpha)
\cdot K_{X'}
\end{Eqn}

\noindent where $\cdot$ denotes intersection product and we identify
the canonical classes with their Poincar\'e duals.  This reflects the
basic ``blow-up'' formula for Gromov invariants in four-manifolds:
the virtual dimension of the space of holomorphic curves in a
four-manifold $W$ in a class $A$ is the same as the dimension of
curves on $W'$ in the class $A+E$, where $E$ is the exceptional
divisor of a blow-up $W' \rightarrow W$.  Indeed the actual invariants
co-incide; the analogue of this for the spaces $X_r(f)$ will be
important later.  To end the section, we shall give an index result.

\begin{Prop} \label{virdim}
Let $\phi(s) = \alpha \in H_2(X';\zz)$ where $s \in \Gamma(F)$.
Then the complex virtual dimension of the space of pseudoholomorphic sections
of $F: X_r(f) \rightarrow \sss^2$ in the homotopy class $s$ is given
by $[\alpha^2 - K_X \cdot \alpha]/2$.
\end{Prop}

\begin{Pf}
By standard arguments,
the virtual dimension is given by the sum of the rank of the vertical
tangent bundle and its first Chern class evaluated on the homology
class of the section.  We can fix a smooth section which passes
through the open dense set $\Sym^r (\Sigma_0 \backslash 
\{ \mathrm{Node} \})$ at each critical fibre, and then we are just
  working with the vertical tangent bundle of a family of symmetric
  products of curves.  It will be helpful to adopt the universal
  viewpoint: via a generic choice of fibrewise metrics,  regard $X'
  \rightarrow \sss^2$ as smoothly embedded inside the total space of
  the universal curve $\mathcal{C}_g \rightarrow \mgbar$.  Indeed, we
  can form the fibre product $\mathcal{Z} = \mathcal{C} \times_{\mu}
  \mathcal{S}^r(\mathcal{C})$ of the universal curve and the universal
  relative Hilbert scheme.  This contains a universal divisor $\Delta$, which
  is just the closure of the obvious fibrewise divisor in $\Sigma
  \times \Sym^r (\Sigma)$ where $\Sigma \in M_g \subset \mgbar$.
  Write $\pi$ for the projection $\mathcal{Z} \rightarrow \mathcal{S}^r
  (\mathcal{C})$ to the universal Hilbert scheme.  Then the vertical
  tangent bundle to $\mu: \mathcal{S}^r (\mathcal{C}) \rightarrow \mgbar$
  is exactly $\pi_* \mathcal{O}_{\Delta} (\Delta)$.  This follows
  from the naturality of the construction of (\cite{ACGH}, pages 171-3).
  From the exact sequence

$$0 \rightarrow \mathcal{O} \rightarrow \mathcal{O}(\Delta)
\rightarrow \mathcal{O}_{\Delta} (\Delta) \rightarrow 0$$

\noindent of bundles on $\mathcal{Z}$, and the associated long exact
sequence in cohomology, it follows that 

$$ch (T^{vt}(\mu)) \ = \ ch (T^{vt} (\Pic)) + ch (\pi_!
\mathcal{O}(\Delta)) - 1.$$

\noindent Here the constant term $1$ comes from $ch(\pi_*
\mathcal{O})$.  If we take the degree zero terms in the above, we get
an equation in the ranks of the bundles:

$$r \ = \ g + (r-g+1) - 1$$

\noindent using Riemann-Roch on each fibre; this checks!  
We have also used the fact that $R^1 \pi_*
\mathcal{O}_{\Delta}(\Delta) = \{ 0 \}$, which holds since skyscraper
sheaves on curves have no higher cohomology.  The term $R^1 \pi_*
\mathcal{O}$ has fibre at $(D \in \Sigma) \in
\mathcal{S}^r(\mathcal{C})$ the cohomology group $H^1 (\Sigma;
\mathcal{O})$, independent of $D$.  This is just the pullback from
$\mgbar$ of the vertical tangent bundle to the Picard bundle. 

\vspace{0.2cm}

\noindent We can apply the
Grothendieck-Riemann-Roch theorem to the term $ch(\pi_!
\mathcal{O}(\Delta))$.  This yields 

$$ch(\pi_! \mathcal{O}(\Delta)) \ = \ \pi_* [ch (\mathcal{O}(\Delta))
\mathrm{Todd}(T^{vt}(\pi))].$$

\noindent We need only work up to degree four on the RHS since we will
evaluate the push-forward on a sphere $[\sss^2]$. The relevant terms
of the Chern character of the line bundle $\mathcal{O}(\Delta)$
are then $1+[\Delta] + [\Delta]^2/2$.  The vertical tangent bundle to
$\pi$ is just (the pullback of)
the vertical tangent bundle to the universal curve.  Hence, the degree
four term on the RHS is given by

\begin{Eqn} \label{equality}
[\Delta]^2/2 - \omega_{\mathcal{C}/M_g} \cdot [\Delta]/2 +
\mathrm{Todd}_{deg(4)} (T^{vt}(\pi)).
\end{Eqn}

\noindent Now evaluate on the sphere arising from the Lefschetz
fibration $X' \rightarrow \sss^2$; by construction, the class $\Delta$
gives $\iota(\alpha)$ whilst the first Chern class of the relative dualising
sheaf gives the canonical class of $X'$ twisted by the canonical
class of the base $\pp^1$:  $c_1 (\omega) = K_{X'} - f^* K_{\pp^1}$.
This twist means that $c_1 (\omega) \cdot
[\Delta]$ is just $\alpha \cdot K_X + 2r$, where $r = \alpha \cdot
[\Fibre]$ inside $X'$.  The third term in (\ref{equality}) cancels
with the term  coming
from $T^{vt} (\Pic)$ above; indeed the relative dualising sheaf is
dual to the vertical tangent bundle to the Picard bundle.  When we put
these pieces together, we find that

$$\mathrm{virdim}_{\cc} (\mathcal{M}_J (s)) \ = \ r + \alpha^2/2 - K_{X}
\cdot \alpha/2 - r $$

\noindent where the first $r$ comes from the rank of the vertical
tangent bundle and the second from the discrepancy between $K_{X'}$ and
$K_{X'} - f^*K_{\pp^1}$.  This gives the required answer.
\end{Pf}


\subsection{Defining the invariants}

We need one last ingredient which was not relevant in the discussions
of \cite{SKDIS}; a method for cutting down the dimensions of moduli
spaces.  Any point $z \in
X$ disjoint from the base-points and critical points of the Lefschetz
pencil gives rise to a smooth divisor $\mathcal{D}(z)$ in the fibre
over $f(z)$ of 
$X_r(f)$.  To define this divisor, note that there is a holomorphic
map

$$\Sym^{r-1}(\Sigma) \ \rightarrow \ \Sym^r(\Sigma); \qquad D \mapsto D+z$$

\noindent whenever $z \in \Sigma$ is a point of a fixed Riemann
surface $\Sigma$.  The image of the map is a smoothly embedded copy of
$\Sym^{r-1}(\Sigma)$ characterised as precisely those points whose
support contains the point $z$.  If we take $\Sigma = f^{-1}(f(z))$
then the image of the map above is exactly $\mathcal{D}(z) \subset
\Sym^r (\Sigma)$, where we suppress the index $r$ for clarity.
Whenever we discuss these divisors $\mathcal{D}(z_i)$, for points
$z_i \in X$, we will assume that the points have been chosen
generically and in particular lie in $X \backslash (\{ p_i\} \cup 
\{ q_j \})$.  Given this
background, we can now define the invariants which play a fundamental role
in this paper.  Recall the maps $i$ and $\iota$ taking $H_2(X)$ to $H_2(X')$.

\begin{Defn} \label{CountSS}
Let $X$ be any symplectic four-manifold and choose a Lefschetz pencil
$f$ on $X$ (always supposed to be of high degree).  Fix $\alpha \in
H_2 (X; \zz)$.   
The \emph{standard surface count} $\mathcal{I}_{(X,f)}(\alpha)$ is
defined as follows: for each $r$, fix some generic $J \in \mathcal{J}_r$.

\begin{enumerate}
\item If $i(\alpha) \not \in \im(\phi)$ then $\mathcal{I}_{(X,f)}(\alpha)
  = 0$; 
\item if $\alpha^2 - K \cdot \alpha < 0$ then
  $\mathcal{I}_{(X,f)}(\alpha) = 0$;
\item If $i(\alpha) = \phi(s)$, with $s \in \Gamma(F_r)$ and
  $[\alpha^2 - K \cdot \alpha]/2 = m \geq 0$, then
  $\mathcal{I}_{(X,f)}(\alpha)$ is the Gromov invariant
  $Gr_{X_r(f)}(s; z_1, \ldots, z_m)$ which counts $J$-holomorphic sections
of $F$ in the class $s$ passing through the fibre-divisors
$\mathcal{D}(z_i), 1 \leq i \leq m$.
\end{enumerate}
\end{Defn}

\noindent Note that it is not \emph{a priori} clear that this is
independent of the choice of Lefschetz pencil; for
applications coming from non-vanishing results that's not actually
important.  The next result explains why the invariant is indeed
well-defined.  Our treatment is formal, since the required Sobolev
machinery is standard; note that we can do without the virtual
fundamental class machinery of
Li and Tian, via a compactness theorem.   The second part of the
Theorem is a ``blow-up'' formula to help with later computations.

\begin{Thm}  \label{blowup}
Let $(X,f)$ be a symplectic Lefschetz pencil (with integral symplectic
  form) and fix a 
  compatible almost complex structure $J$ on $X$.  Use this to define
  the associated fibre bundles as above, and equip these with
  almost complex and symplectic structures.  Fix $\alpha \in H_2(X; \zz)$.

\begin{itemize}
\item The invariant $\mathcal{I}(\alpha)$ (when not
  defined to be zero) can be computed as the signed count of the
  points of a compact zero-dimensional moduli space in which each point has a
  uniquely attached sign $\pm 1$. 

\item $i(\alpha) \in \im(\phi) \Leftrightarrow
  \iota(\alpha) \in \im(\phi)$ and in obvious notation $\mathcal{I}(\alpha) =
  \mathcal{I}(\iota(\alpha))$. 
\end{itemize}
\end{Thm}

\begin{Pf}
Most of this is implied by the arguments of
\cite{SKDIS}. Let us review the
principal parts.  

\vspace{0.2cm}

\noindent $\bullet$ For generic $J$ on $X_r(f)$, compatible with
the fibration (i.e. $J \in \mathcal{J}$, inducing the given integrable
structures on the 
fibres), the moduli space of holomorphic sections in a given homotopy
class $s$ will be a smooth open manifold of real dimension 

$$d(s) = C_{\phi(s)}\cdot C_{\phi(s)} - K_{X'} \cdot C_{\phi(s)}.$$

\noindent According to general theory, we can compactify the space by
adding two kinds of element.
First, the cusp sections:  that
is, holomorphic sections which may have a number of bubbles,
necessarily lying in fibres of $F: X_r(f) \rightarrow \sss^2$. Second,
we add curves
which are not actually \emph{sections} in the sense that we may add
curves passing through the critical values of $F$.  We need to know
that the points we add do not have excess dimension.  For the cusp
curves, argue as follows.  Any bubble in any fibre projects to a sphere in
the Picard fibration $P_r(f)$.  This either lies in the smooth part of
a fibre or lifts
to the normalisation of a singular fibre.  Arguing as in 
Lemma (8.11) of \cite{SKDIS} - which applies unchanged to the case of
general $r$ and not just $r=2g-2$ - the latter situation cannot occur
and hence any bubble represents a
multiple of the generator of a projective space fibre for some linear
system.  If we split off $n$ bubbles, the remaining section component
$s'$ maps under $\phi$ to the homology class $C_{\phi(s)} -
n[\Fibre]$.  Provided at least the degree of the pencil is large, so
$C \cdot [\Fibre] \gg 0$, the virtual dimension for holomorphic
sections in this new class will be very negative.  Then, by
the assumed regularity of $J$, when the degree
of the pencil is large enough the total
dimension of the space of cusp curves is negative and hence
\emph{there are no bubbles in any moduli spaces}. 

\vspace{0.2cm}

\noindent To see that the curves passing through the
critical values cannot be of excess dimension, the argument is even
stronger.  Namely, we claim that
\emph{no curve passing through the critical loci can arise as the
  limit of a sequence of smooth sections} for \emph{any} $J$.  For
although the limit
curve is (in the appropriate Sobolev setting) only a $L^{2,2}$-map,
the section component $s$ is a holomorphic curve for a smooth almost
complex structure on an almost complex manifold and hence is smooth by
elliptic regularity.  The composite map $F \circ s: \sss^2 \rightarrow
\sss^2$ is then holomorphic and degree one away from the
isolated point of intersection with the normal crossing divisor in the
singular fibre.  Hence the map is a degree one diffeomorphism
everywhere, and hence $s$
is indeed a section; so it cannot pass through any point of the total
space where $dF=0$.  We have already ruled out bubble components, and
hence \emph{there are no curves passing through the singular loci}.

\vspace{0.2cm}

\noindent We deduce that the moduli space of holomorphic sections 
carries a fundamental class of the correct (virtual) dimension.
Global orientability of the moduli space follows since the
$\cdbar$-operator of a non-integrable $J$ is just a zeroth order
perturbation of the usual $\cdbar$-operator.  The standard surface
count, as we have defined it, is given by cutting down the moduli
space to be zero-dimensional.  Note that cutting down the moduli
spaces introduces no 
analytic problems since we are only dealing with spheres, in which
case evaluation maps are always submersions by a result of
\cite{McD-S:Jhol}.  At this stage, orientability is exactly
the assignment of a sign to each point, giving the first part of the
above theorem.  That we really have an invariant, independent of
the choice of almost complex structure, follows from the usual cobordism
argument for one-parameter families together with the observation
(\ref{connectedparameters}). 

\vspace{0.2cm}

\noindent $\bullet$ To see that it is equivalent to count holomorphic sections
yielding surfaces in a class $i(\alpha)$ or $\iota(\alpha)$ in
$H_2(X';\zz)$, we use the analogue of the divisors $\mathcal{D}(z)$
coming from the exceptional sections.  That is, notice an exceptional
curve $E$ defines a point $e_t \in \Sigma_t$ for each $t \in \sss^2$
and hence a smooth submanifold $X^E_{r-1}(f) \hookrightarrow X_r(f)$.
As in \cite{SKDIS}, there are almost complex structures on $X_r(f)$
which give rise to smooth almost complex structures on each
$X^E_{r-1}(f)$ and for which the natural inclusions are holomorphic.
Choose $J$ on $X_r(f)$ regular amongst almost complex structures in
$\mathcal{J}$ with this property; then the compactified space of $J$-sections
yielding curves in the class $\iota(\alpha)$ will be a smooth compact
manifold where none of the curves contain bubbles.  We claim that in
fact all of these sections lie inside the intersections of the images
of all the $X^E_{r-1}(f)$, and that this defines an isomorphism of
moduli spaces $\mathcal{M}_J (\iota(\alpha)) \sim \mathcal{M}_J
(i(\alpha))$.  Suppose otherwise; then some element $s \in
\mathcal{M}_J (\iota(\alpha))$ does not lie inside $X^E_{r-1}(f)$, and
hence meets it with locally positive intersection.  Hence the cycle
$C_s$ in $X'$ meets $E$ with locally positive intersection; but the
algebraic intersection number is just $E \cdot [C_s] = E \cdot [\alpha
+ \sum E_i] = -1$, and this yields a contradiction.  

\vspace{0.2cm}

\noindent So all the sections lie inside all the loci defined by
exceptional sections, and hence give rise to canonical holomorphic
sections of the intersection $\cap_E X^E_{r-1}(f) \sim X_{r-N(E)}(f)$,
where $N(E)$ is the number of exceptional curves.  But then
$i(\alpha)\cdot[\Fibre] = \iota(\alpha)\cdot[\Fibre]-N(E)$, and it is
easy to check this yields the required isomorphism of moduli spaces.
The virtual dimensions for the classes co-incide, so we cut down with
a fixed set of $\mathcal{D}(z_i)$ to obtain isomorphic compact
zero-dimensional moduli spaces.  The reader can easily check that the
signs associated to points agree; the result follows.
\end{Pf}

\noindent Let us draw attention to a part of the above, a ``fibred
monotonicity'' property:

\begin{Cor}
For a pencil $f$ of high degree on $X$, and any generic almost complex
structure $J \in \mathcal{J}$ on $F:X_r(f) \rightarrow \sss^2$, all
moduli spaces 
of smooth holomorphic sections of $F$ are already compact.
\end{Cor}

\noindent The above arguments show that whenever we
have a non-zero standard surface count for a class $\iota(\alpha)$ we
can find symplectic surfaces in $X'$ disjoint from the $E_i$ and hence
push them down to surfaces in $X$.  The importance of this is that it
allows us to stabilise the intersection number $r$ to be large, and
hence take advantage of the geometry of the Abel-Jacobi map.  On the
other hand, the blow-up identity is also important: for instance,
$\iota(\alpha)$ and $i(K_X - \alpha)$ meet the fibres of a Lefschetz
pencil respectively $r$ and $2g-2-r$ times, and this is where Serre
duality will enter.  Indeed, the explicit
computation of the standard surface count for the canonical class in
\cite{SKDIS} involved dealing with a class of almost complex structures
- those \emph{compatible with the zero-sections} - which were tailored
to the Abel-Jacobi map, and it is a strict
generalisation of this latter class of structures which will underlie
Theorem (\ref{Serreduality}).


\section{Standard surfaces and holomorphic curves}

In the previous section we defined new invariants of Lefschetz
pencils, but did not relate them to symplectic submanifolds in $X$.
The point is 
that, for appropriate complex structures, the sections of $X_r(f)$
provided by the non-triviality of a standard surface count yield
standard symplectic surfaces.  In \cite{SKDIS} we reached this
conclusion via an intermediate stage:  we constructed unions of such
surfaces with positive local intersections and then applied a smoothing
lemma.  This has the advantage of simplicity, but the disadvantage of
remaining bound within the realm of \emph{integral} symplectic
manifolds:  we always assumed the ray $\rr_+ \langle \omega_X \rangle$
defined a point of $\pp H^2 (X; \qq) \subset \pp H^2(X;\rr)$.  To
avoid this assumption, here we will build holomorphic curves in $X$
from non-triviality of the $\mathcal{I}$-invariants.  Additionally,
we will sketch a proof of the result 
``$\mathcal{I}=Gr \ (\mathrm{mod} \ 2)$'', at least in the absence of
square zero tori.  This result is intended
to be motivational, and to make the present paper coherent,  more
than to make serious headway on a full proof of (\ref{IequalsGr}).  


\subsection{Almost complex diagonals}

To obtain symplectic surfaces, we need to use special almost
complex structures on $X_r(f)$, and for this we need to introduce the
\emph{diagonal strata} (cf. \cite{SKDIS}, Section 6).   For each
partition $\pi$ of the form $r = \sum a_i n_i$ and each smooth fibre,
we have a map  

$$\Sym^{n_1}(\Sigma) \times \cdots \times \Sym^{n_s}(\Sigma) \
\longrightarrow \ \Sym^r (\Sigma); \qquad (D_1, \ldots, D_s) \mapsto
\sum a_i D_i.$$

\noindent We induce a smooth map of fibre bundles 

$$Y_{\pi} \ = \ X_{n_1}(f) \times_f \cdots \times_f X_{n_s}(f) \ \rightarrow
X_r(f)$$

\noindent which is finite and generically a homeomorphism onto its
image.  An almost complex structure $J$ on $X_r(f)$ is
\emph{compatible with the strata} if there are almost complex
structures $j_{\pi}$ on $Y_{\pi}$, for every partition $\pi$, making
the above maps $(j_{\pi}, J)$-holomorphic.  That these exist
follows from a local computation with restricted charts:  the point is
that almost complex structures on each of $X_r(f), Y_{\pi}$ arise by
patching local canonically defined structures via partitions of
unity.  Since a smooth partition of unity on $X_r(f)$ defines smooth
partitions of unity on all the $Y_{\pi}$, pulling back under the
canonical smooth maps above, we deduce existence.  In fact, we obtain the
following, given as
Propositions (6.3) and (7.4) in \cite{SKDIS}:

\begin{Lem}
Let $\mathcal{J}_{\Delta} \subset \mathcal{J}$ denote the class of
smooth almost complex
structures on $X_r(f)$ which are compatible with the strata.  Then
this is non-empty, and there is an open dense set $\mathcal{U}$ in
$\mathcal{J}_{\Delta}$ with the following property.  If $J \in
\mathcal{U}$ then all
moduli spaces of smooth holomorphic sections of all fibrations $Y_{\pi}$,
including $X_r(f)$ itself, have the expected dimension.  Moreover, a
dense set of points of each moduli space corresponds to sections
transverse to all diagonal strata in which they are not contained.
\end{Lem}

\noindent The regularity again follows standard lines:  locally,
perturbations of almost complex structures compatible with the strata
generate the whole tangent space to the space of sections (i.e. all
vector fields tangent to the fibres along the image $\sss^2$).  Note
that we can immediately assert:

\begin{Prop}
If $s: \sss^2 \rightarrow X_r(f)$ is a section which meets the
diagonal strata transversely at embedded points and with locally
positive intersections, then the cycle $C_s \subset X'$ is a (not
necessarily connected) smooth standard surface.
\end{Prop}

\noindent This is just because transverse intersections with the
top stratum of the diagonals (and no other intersections) give
tangency points of $C_s$ to the fibres of $f$, whereas singularities
of $C_s$ arise from non-transverse intersections with the diagonals;
for instance, nodes of $C_s$ arise from tangencies to the diagonals.
Before continuing, we need further  discussion of almost complex
structures on fibre bundles and their associated symmetric product
bundles.

\vspace{0.2cm}

\noindent For any fibration of manifolds $\pi:Z \rightarrow B$ we can form
the fibre product $Z \times_{\pi} \cdots \times_{\pi} Z \rightarrow B$
whose fibre over $b \in B$ is just $\pi^{-1}(b) \times \cdots \times
\pi^{-1}(b)$.  That this is smooth follows immediately from the
surjectivity of $d(\pi \times \cdots \times \pi)$ viewed as a map $Z
\times \cdots \times Z \rightarrow B \times \cdots \times B$.  If the
original manifolds $Z$ and $B$ carry  almost complex structures $J,j$ for
which $\pi$ is a holomorphic projection, then the product complex
structure $J \times \cdots \times J$ induces an almost complex
structure on the total space of the fibre product.  Now if we remove
all the diagonals from the fibre product, then there is a free and
holomorphic action of the symmetric group $S_r$, and the quotient
inherits a natural smooth almost complex structure: in other words,
for any pseudoholomorphic fibration, the space $Z_r(\pi) \backslash
\{\mathrm{Diagonals}\}$ carries a natural smooth almost complex
structure $\mathbb{J}$.

\vspace{0.2cm}

\noindent In general we cannot say any more that this;  the total
space of the relative symmetric product $Z_r(\pi)$ will not carry the
structure of a smooth manifold, even, and we cannot make sense of the
tangent space at the points of the diagonal.  If however the fibres of
$\pi$ are two-real-dimensional surfaces, then using families of
restricted charts, we find that the fixed $J$ on $Z$ does
induce on $Z_r(\pi)$ the structure of a smooth manifold.  In this
case, we can ask about the induced almost complex
structure $\mathbb{J}$ in a neighbourhood of the diagonals; from the
topological set-up it is clear that it extends continuously, but
higher regularity is not obvious.  Let $J$ be given in local complex
co-ordinates $z,w$ on $X'$ by a matrix 

$$\left( \begin{array}{cc} i & \mu \\ 0 & i \end{array} \right)$$

\noindent for a complex anti-linear homomorphism $\mu: f^* T_b \sss^2
\rightarrow T_p X'$, so $\mu \cdot i + i \cdot \mu = 0$.  Here $z$ is
the fibre co-ordinate and $w$ a co-ordinate in the base.  Taking
product charts, reduce to the case of a point 
$p+\cdots+p$ in the small diagonal of the symmetric product.  There are
induced co-ordinates $\sigma_1, \ldots, \sigma_r, w$ for $\sigma_i$
the $i$-th elementary symmetric function of $r$ copies $z_i$, $1 \leq
i \leq r$, of the local co-ordinate $z$ near $p$.  Then $\mathbb{J}$ is given
near this point as 

$$\left( \begin{array}{cc} \cdots & \mu_1 \\ \Sym^r(i) & \vdots \\ \cdots &
    \mu_r \\ 0 & i \end{array} \right).$$

\noindent The functions $\mu_i$ are obtained from smooth functions of
$\mu$ which are invariant under the action of $S_r$ on the fibre
product.  The elementary symmetric functions $\sigma_i$ are polynomial
in given co-ordinate functions $x_i, y_i$ upstairs, and hence the smooth
function $\mu$ of the real co-ordinates is \emph{a priori} just a
function of the 
fractional powers $\sigma^{1/i}$ (and their complex conjugates) for $i
\leq r$.   Nonetheless, some regularity does persist.
The following was observed independently by the
author and, in a mildly different context, by Siebert and Tian
\cite{Sieb:weier}.  

\begin{Lem} \label{Holder}
The extension of the almost complex structure $\mathbb{J}$ on
the relative symmetric product $Z_r(\pi)$ from the
complement of the diagonals to the total space is H\"older
continuous, of H\"older exponent $C/r$ for some constant $C$.  Along
the top open stratum where at most two points co-incide the extension
is Lipschitz.
\end{Lem}

\begin{Pf}
The key computation, due to Barlet, is the following.  Be given a
holomorphic branched  
covering $D \rightarrow D'$ of  complex domains of maximal ramification order
$n$ and with branch locus having normal crossing singularities.  Let
$\phi$ be a smooth function on $D$ and let $\overline{\phi}$
be the ``trace'' of $\phi$ on $D'$, that is $\overline{\phi}(z) =
\int_{\phi^{-1}(z)} \phi$.  Then $\overline{\phi}$ is H\"older
continuous of exponent $2/n$.  As a special case, if a group $G$ acts
holomorphically 
on a complex manifold $D$ and $D / G$ is smooth, then a smooth
$G$-invariant function $\phi$ on $D$ defines a H\"older continuous
function on the quotient.  Since any holomorphic
branched covering can be resolved to have normal crossing branch
locus, it follows that traces of smooth functions are always H\"older
continuous of exponent depending only on the ramification.  In our
situation this applies to the branched covering from the fibre product
to the symmetric product, ramified along the diagonals.  If only two
points co-incide, then the ramification locus is already smooth: there
is no need to resolve, and the extension is H\"older continuous for
every $\alpha < 1$, in other words Lipschitz continuous.  The local
computation - and much besides - is given in Barlet's paper
\cite{Barlet} (this reference is due to Siebert and Tian). 
\end{Pf}

\noindent Notice that the above implies a \emph{tautological
correspondence} between smooth $J$-holomorphic curves in $X$ and
H\"older $\mathbb{J}$-holomorphic sections of $X_r(f)$ (in
the case where the 
sections lie inside a diagonal stratum, one can check this by pulling
back to the minimal such stratum, where the complex structure
$\mathbb{J}$ is now generically smooth).  Such a correspondence
presumably lies at the heart of (\ref{IequalsGr}).  


\subsection{Relation to the Gromov invariant}

Basic theory of, and crucially the compactness theorem for,
pseudoholomorphic curves has been proven with rather weak assumptions
on the regularity of the almost complex structure: the best reference
is \cite{AudLaf:eds}, in particular the article by J.C. Sikorav.

\begin{Prop}
Let $X$ be a symplectic four-manifold and fix a taming almost
complex structure $J$ on $X$.  Suppose for every Lefschetz
pencil $f$ of sufficiently high degree the standard surface count
$\mathcal{I}_{(X,f)}(\alpha) \neq 0$.  Then there is a $J$-holomorphic
curve in $X$ in the homology class $\alpha$.
\end{Prop}

\begin{Pf}
We can choose a sequence $J_n$ of almost complex structures on $X$
which converge in (say) $C^2$-norm to $J$ and such that for each
$n$, there is a Lefschetz pencil of $J_n$-holomorphic curves on
$X$.  (The degrees of these pencils may have to increase with $n$,
hence the wording in the hypothesis of the Proposition.)
This follows from the main theorem of \cite{Donaldson:submflds}; for any fixed
$\varepsilon > 0$ we know $X$ admits Lefschetz pencils of surfaces
whose tangent spaces deviate from integrability by at most
$\varepsilon$ measured in any given $C^k$-norm.  Hence it will be
enough to show that for an almost complex structure $j$ on $X$ for which a
holomorphic Lefschetz pencil (of high enough degree) exists, there are
$j$-holomorphic curves in the class $\alpha$.  Then we can finish the
proof using Gromov compactness for a sequence of $J_n$-holomorphic curves.

\vspace{0.2cm}

\noindent Now fix such a $j$ on $X$ and form the fibre bundle $X_r(f)
\rightarrow \sss^2$, where $r = \alpha\cdot[\Fibre]$ in the usual
way.  The almost complex structure $j$ on $X$ induces a
canonical smooth almost complex structure on $\mathcal{X}^* =
X^*_r(f)\backslash\Delta$, 
that is on the complement of the diagonals in $X^*_r(f) =
X_r(f) \backslash \{F^{-1}(\Crit(f))\}$.  To see this,
note that $X' \rightarrow \sss^2$ is a smooth fibre bundle away from
the critical fibres; now the above discussion 
yields such a canonical $\mathbb{j}$ on  $\mathcal{X}^*$, at least
H\"older continuous at the diagonals. Near the critical fibres, the
original fibration was actually holomorphic, and the relative Hilbert
scheme carries a smooth integrable complex structure. The data
patches, by naturality, and so we obtain a H\"older continuous almost
complex structure $\mathbb{j}$ on $X_r(f)$ which is smooth on a
dense set.

\vspace{0.2cm}

\noindent The smooth sections of any vector bundle are always dense in
the H\"older sections.  It follows that we can
choose a sequence of smooth almost complex structures
$\mathbb{j}_n$ on $X_r(f)$ which converge to the canonical structure,
in $C^{0,\alpha}$-norm say.  Since
the structures on $X_r(f)$ compatible with the strata form a dense
Baire set, we can assume that each $\mathbb{j}_n$ lies in
$\mathcal{J}_{\Delta}$.  For each $n$ we then have a section $s_n$ of
$X_r(f)$ in the homotopy class $\phi(\alpha)$, and the results of
\cite{SKDIS}, as in the previous section, assert that this defines a
positive symplectic divisor
$C'_n$ in $X'$ which contains the exceptional curves and descends to a
symplectic divisor $C_n$ in $X$ in the class $\alpha$.  The symplectic
condition controls the genera and area of all the surfaces uniformly.

\vspace{0.2cm}

\noindent For each surface $C_n$ we can find an almost complex
structure $j_n$ on $X$ for which $C_n$ is $j_n$-holomorphic.  Since we
have convergence of the $\mathbb{j}_n$ in
$C^{0,\alpha}$ on $X_r(f)$, we can choose the $j_n$ so as to converge
in $C^{0,\alpha}$-norm to the given 
almost complex structure $j$.  In this case, we can complete the proof
using the Gromov compactness theorem; this holds for sequences of
almost complex structures converging only in $C^{0,\alpha}$ by the
results of Sikorav and Pansu in \cite{AudLaf:eds}.  The
compactness will
yield a $j$-holomorphic curve $C$ in the class $\alpha$ (for the proof of
Pansu shows that no area is lost in the limit), and since the almost
complex structure $j$ is smooth, elliptic regularity asserts that the
curve $C$ is the image of a smooth map. 
\end{Pf}

\noindent  Let's add in the following well-known regularity result, due to
Ruan (and developed by Taubes):

\begin{Prop}[Ruan]
Let $X$ be a symplectic four-manifold and $\alpha \in H_2 (X;\zz)$.
There is an open dense set
$\mathcal{J}_{reg}$ in 
the space of compatible almost complex structures on $X$ for which the
following hold for $j \in \mathcal{J}_{reg}$:

\begin{itemize}
\item if $[\alpha^2 - K_X \cdot \alpha]/2 = d \geq 0$ then the space of
  $j$-holomorphic curves
  representing $\alpha$ and passing through $d$ generic points of $X$
  is finite; 
\item each such holomorphic curve is embedded or an unramified cover
  of a square zero torus.
\end{itemize}
\end{Prop}

\noindent An immediate consequence of this, and the preceding
existence result, is the following

\begin{Cor}
Let $(X, \omega)$ be a symplectic four-manifold.  Let $\omega'$ be a
rational  perturbation of $\omega$ and construct Lefschetz pencils $f$
on $X$ with fibres dual to $k[\omega']/2\pi$. 
If $\mathcal{I}_{(X,\omega',k)}(\alpha) \neq 0$ for all $k \gg 0$ then
$\alpha$ may be
represented by embedded $\omega$-symplectic surfaces in $X$.
\end{Cor}

\noindent Since there is no rationality hypothesis on
the original symplectic form on $X$, this statement - even for
the canonical class - is stronger than that obtained in \cite{SKDIS}.

\vspace{0.2cm}

\noindent  Unfortunately, H\"older regularity of the almost complex
structure $\mathbb{J}$ does not seem strong enough to directly apply
the usual implicit 
function theorem.  Hence, it seems non-trivial to prove that the map
from holomorphic sections of $F$ to holomorphic curves in $X$ outlined
above is \emph{onto}.  This would prevent any naive
comparison of the $\mathcal{I}$ and $Gr$ invariants.  However, in
principle we can get around this, using the fact that along the
top stratum of the diagonal we actually have Lipschitz control.  First
we need to strengthen Ruan's regularity result to a fibred situation.
An almost complex structure on a four-manifold is \emph{compatible}
with a Lefschetz pencil $f$ if locally near each base-point there is a
$j$-holomorphic diffeomorphism to a family of complex discs through
the origin in $\cc^2$ and if all the
fibres of $f$ are $j$-holomorphic curves in $X$.  Considering the
sphere in moduli space defined by the pencil, for instance, we find 
that such $j$ always exist and form a connected subspace
$\mathcal{J}_f$ of the space of all $\omega$-taming forms.

\begin{Lem} \label{regularj}
Let $X$ be a symplectic four-manifold and fix $\alpha \in
H_2(X;\zz)$.  For Lefschetz pencils $f$ of sufficiently
high degree on $X$ and for almost complex structures $j$ on $X$
generic amongst those compatible with $f$, the moduli spaces of
$j$-holomorphic curves in the class $\alpha$ are smooth and of the
correct dimension. 
\end{Lem}

\begin{Pf}  Fix some pencil $f$ of symplectic surfaces of area
  greater than $\alpha \cdot [\omega]$.  If $j$ varied in the class of
  all taming almost complex
  structures on $X$ this would be exactly Ruan's theorem, but we are
  restricting to almost complex structures for which the fibres of the
  pencil $f$ are holomorphic.  In particular, such a $j$ cannot be
  regular for curves in the class of the fibre.  Nonetheless, by
  positivity of intersections of holomorphic curves, any
  $j$-curve $C$ in the class $\alpha$ will be some multisection of the
  fibration (by consideration of area no such curve can contain
  any fibre components).  In local co-ordinates near some point $p \in
  X$ disjoint from the base locus of $f$, all perturbations of $j$ in
  horizontal 
  directions (i.e. in the term $\mu$ in the description before
  \ref{Holder}) yield almost complex structures which are compatible
  with $f$.  Hence the class of compatible $j$ is large enough for
  deformations $(C',j')$ of $(C,j)$ (through curves holomorphic for
  some $j' \in \mathcal{J}_f$) 
  to generate the entire tangent space to $X$, at
  least over an open subset where the given curve is a multisection and
  not tangent to any fibre.  By an argument with Aronszajn's Theorem,
  this is enough to 
  achieve transversality, as in \cite{McD-S:Jhol} or the related
  result (\cite{SKDIS} Lemma 7.5). 
\end{Pf}

\noindent For generic almost complex structures on a four-manifold
with $b_+ > 1$,  
the Gromov invariant used by Taubes counts all curves that are not
square zero tori with weight $\pm 1$, whilst the Gromov invariant counting
somewhere injective spheres  used to define
the standard surface count - cf. \cite{McD-S:Jhol} for instance  -
again counts each point of the moduli
space with weight $\pm 1$.  Hence a bijection of moduli spaces will
give (\ref{IequalsGr}) modulo two, at least for $X$ containing no
embedded symplectic tori of square zero and for 
zero-dimensional moduli spaces.  (In 
fact, when $b_+ > 2$, we know that the only non-trivial standard
surface counts do occur in 
zero-dimensional moduli spaces; we will deduce this from
(\ref{Serreduality}) later in the paper, and the same property is
standard for the Gromov invariants.  But it is no harder to
establish the bijection 
even if we first have to cut down dimensions.)   The previous lemma
shows, in essence, that there are $f$-compatible $j$ which
are sufficiently regular for the computation of the Gromov invariant
of the class $\alpha$.  By
perturbing $j$ amongst $f$-compatible structures again if necessary, we see:

\begin{Cor} \label{transverse}
Suppose $X$ contains no symplectic square zero tori.  Fix $\alpha \in
H_2 (X;\zz)$ and a Lefschetz pencil $f$ on $X$ as
in (\ref{regularj}).  For generic $j$ compatible with $f$, the
invariant $Gr_X(\alpha)$ is given by a signed count of finitely many
embedded 
holomorphic curves, each of which is transverse to the fibres of $f$
away from finitely many tangencies all at smooth points.
\end{Cor}

\noindent For this $j$, the associated $\mathbb{J}$ on the relative
Hilbert scheme has - by the tautological correspondence - a moduli
space of sections in the class $\phi(\alpha)$ comprising finitely many
points, each of which is a section which meets the diagonal strata
transversely (and hence only in the top stratum).  We need a key
technical result, due to Siebert and Tian (\cite{Sieb:weier}, Theorem
II): we have phrased their result in terminology appropriate to our
situation.  Fix a restricted chart $D \times D \rightarrow X$ and
work, by pulling back, locally in the domain.  Hence we start with a smooth
family $\pi: D \times D \rightarrow D$ of complex discs over the disc,
with a smoothly varying family of
holomorphic structures on the fibres.  As usual, there is a fibration
$\cc_r(\pi) \rightarrow D$ with fibres the $r$-th symmetric products
of fibres of $\pi$, and this carries an induced almost complex
structure $\mathbb{J}$ with regularity as in (\ref{Holder}).  Fix some
finite $p>2$.

\begin{Thm}[Siebert, Tian] 
If $C$ is a $\mathbb{J}$-holomorphic section of $\cc_r(\pi)
\rightarrow D$ with finitely many
transverse intersections with the diagonals, then the moduli space
of $\mathbb{J}$-holomorphic sections near $C$ is naturally a smooth
Banach manifold, modelled on the Banach space $L^{1,p}_{hol}(D;\cc^r)$.
\end{Thm}

\noindent In other words, we can parametrise local pseudoholomorphic
curves provided the almost complex structure is at worst Lipschitz and
smooth off a finite set.  (Indeed one can relax the last condition.)
In \cite{AudLaf:eds} Sikorav observes that the $\cdbar$-equation can
be linearised at any H\"older complex structure, and that
solutions $f$ of $\cdbar_j (f) = 0$ - with $j$ H\"older of exponent
$\alpha$  - will
necessarily be of class $C^{1,\alpha}$ and hence differentiable
everywhere.  The key in the result above is to show that, in the
Lipschitz setting, an approximate
right-inverse to the linearisation (provided as usual by the Cauchy
integral operator) satisfies strong enough bounds to
apply a contraction mapping and parametrise local $\mathbb{J}$-curves
by perturbations of the given curve $C$.  

\vspace{0.2cm}

\noindent Now suppose we have a $\mathbb{J}$-holomorphic section $s:
\sss^2 \rightarrow X_r(f)$ 
which satisfies the transversality properties coming from
(\ref{transverse}).  Combining the local perturbations above,
near its intersections with the diagonals, with
the usual perturbation arguments for holomorphic curves in smooth
almost complex manifolds, we can see that globally over the base of the
Lefschetz pencil the linearisation of the
$\cdbar_{\mathbb{J}}$-operator will be surjective at $s$.  For the argument
of \cite{McD-S:Jhol} (page 35, cf. the footnote) works here even
without \emph{a priori} $L^{2,2}$-regularity on the maps and
Aronszajn's Theorem.  That is, an element
in the kernel of the adjoint of the linearisation must be orthogonal
to all vertical vector fields (since a given
section of $X_r(f)$ is never tangent to the fibres), and hence
trivial.  In this context, on the open subset $X^0_r(f)$ of $X_r(f)$
given by removing a neighbourhood of the
lower strata of the diagonals, we 
can use the implicit function theorem for the
$\cdbar$-equation.  This is viewed as a section of the bundle over the space
$C^{1,\alpha}(\sss^2; X^0_r(f)) \times \mathcal{J}^{Lip}$ with fibre 
the one-form valued sections $C^{0, \alpha}(\sss^2, u^* (TX^0_r(f))^{vt}
\otimes K_{\sss^2})$.  The implicit function theorem in particular
shows that the moduli space of holomorphic sections of $X_r^0(f)$ is
independent, up to diffeomorphism, of sufficiently small Lipschitz
perturbations of the almost complex structure.  On the other hand, we
know that in the space of all H\"older almost complex structures on
$X_r(f)$ the smooth structures are dense, and that in this dense
subspace the structures for which all points of zero-dimensional
moduli spaces represent sections having only transverse intersections
with the diagonals are also dense.  Hence, approximating $\mathbb{J}$
by smooth structures $J_n$ (convergent in $C^{0, C/r}$) which satisfy
this genericity, which
\emph{a priori} admit no holomorphic sections lying inside the
diagonals, and which Lipschitz approximate $\mathbb{J}$ elsewhere,
yields the following. 

\begin{Prop}
Fix a generic $j$ on $X$ compatible with a sufficiently high degree
pencil $f$.  The moduli space of $\mathbb{J}$-holomorphic sections of
$X_r(f)$ in a homotopy class $\phi(\alpha)$ is a smooth Banach
manifold.  Moreover, this manifold is diffeomorphic to the moduli
space of $J_n$-holomorphic sections for any sufficiently large $n$,
where $(J_n)_{n \in \nn}$ is a $C^{0, C/r}$-approximating sequence of smooth
almost complex structures on $X_r(f)$ as above. 
\end{Prop}

\noindent Here, by cutting down with point conditions in $X$ or with
the $\mathcal{D}(z)$ in $X_r(f)$, we are implicitly working only with
zero-dimensional 
moduli spaces.  In this setting, the smoothness of the moduli space
amounts to the regularity of the
isolated points (they carry no obstruction information).  The partial
result (\ref{IequalsGrmodtwo}) is almost a 
consequence; by our general theory, the ``mod two'' Gromov invariant
on $X$ is the  parity of the number of
points in the moduli space of holomorphic sections of $X_r(f)$ for
the canonical H\"older almost complex structure.  This is now identified
with the sections for a nearby smooth almost complex structure.
The signed count of the second of these finite sets is by
definition the standard surface count of the class $\alpha$.  It's not
immediate from the discussion so far
that any fixed pencil $f$ is of sufficiently high degree to obtain
equality of $\mathcal{I}$ and $Gr$ modulo two for all classes
simultaneously.  But in fact each of the invariants can be non-zero for at
most finitely many classes, and
so if we choose $f$ as in (\ref{regularj}) varying $\alpha$ over
a suitable finite set we obtain the theorem.

\begin{Rmk} Recall that for any fibre bundle $Z \rightarrow B$ which has
holomorphic total space and projection, the fibrewise symmetric
product is naturally almost complex away from the diagonals.
Suppose again we have two-real-dimensional fibres and a global smooth
structure does exist.  Fix a connexion on
$f:Z \rightarrow B$, so at each point $p$ we have a splitting

$$T_p Z \ = \ T_p F_{f(p)} \oplus \langle Hor_p \cong f^* T_{f(p)}B \rangle.$$

\noindent If $Z_r(f) \rightarrow B$ is smooth, then it carries a
distinguished connexion: at a tuple $p_1 + \cdots + p_r$ of points in
one fibre $F_b$ we have the
tangent space to the fibre and a splitting given by the diagonal

$$T_b B \cong \Delta \subset Hor_{p_1} \oplus \cdots \oplus Hor_{p_r}.$$

\noindent Now the connexion induces an almost complex structure
$\mathbb{j}$ on $Z_r(f)$ by taking the product structure in charts, so
the holomorphic sections are just the flat sections;  if the original
almost complex structure $j$ on $Z$ is itself split, then this serves
to extend the canonical structure on $\mathcal{Z}^*$, as described
above, over the diagonals.  (More simply, if $\mu = 0$ before
(\ref{Holder}) then all the $\mu_i = 0$, and these are smooth!)

\vspace{0.2cm}

\noindent For a Lefschetz fibration, there is a distinguished
connexion defined by the symplectic form, but this does not extend
over the critical points;  equivalently, the flat almost complex structure
$j$ on $X^* \rightarrow \pp^1$ defined on the complement of the
critical fibres does not extend to a smooth structure on $X'$.
However, for smooth fibred four-manifolds there is no problem, and one
can use the above point of view to prove (\ref{IequalsGrmodtwo}) in
this special case; the almost complex structures are smooth but
we do not have regularity, so instead we need to analyse the
obstruction theory.
\end{Rmk}

\noindent Proving the equivalence ``$\mathcal{I} = Gr$''
seems particularly appealing considering the fact
that, to date, there is no (non-homotopy) invariant of a general
symplectic manifold which is defined from one Lefschetz pencil and
proved independent of the choice of pencil.  Note also that this would
give a version of the Gromov-Taubes invariant of the four-manifold in
which one only had to count \emph{spheres}, for which the analysis is
often simpler than for counting higher genus curves.


\subsection{Heuristic interlude}

In this subsection we give three digressions to related ideas.  First,
with future higher-dimensional applications in mind, we sketch how to
obtain holomorphic curves from geometric measure theory.  Second, we
outline one approach to proving the standard surface counts are
independent of the choice of Lefschetz pencil.  Lastly, we discuss
some ideas relating the conjectural equivalence of these invariants
and the Seiberg-Witten invariants.

\vspace{0.3cm}

\noindent $\bullet$ For any symplectic manifold,
an almost complex structure $J$ is \emph{strictly compatible} if it
tames the symplectic form and if $\omega(Ju,Jv) = \omega(u,v)$ for any
tangent vectors $u,v$.  The space of these is contractible.
For constructions of
surfaces in fibred six-manifolds (where all the fibres are Hilbert
schemes and H\"older regularity is less clear) another strategy may be
important, so we sketch the argument;  it also appears in Taubes
\cite{Taubes:SWtoGr}.  

\vspace{0.2cm}

\noindent Recall that, by approximating the canonical almost complex
structure on the relative Hilbert scheme by smooth structures
compatible with the strata, we had a
sequence of symplectic surfaces $C_n$ of uniformly bounded genus and
fixed homology class. There is therefore a
fixed upper bound $K$, independent of $n$, on the number of points of
$C_n$ which lie in
one of the critical fibres of $f$ or which are tangent to one of the
fibres.  Reparametrising, we can assume we are in the following
situation:  we have a fixed surface $\Sigma$ and a sequence $h_n:
\Sigma \rightarrow X$ of symplectic immersions of $\Sigma$ in a fixed
homology class $\alpha$.  Moreover there is a finite set of $K$ points
on $\Sigma$ such that the surfaces $h_n (\Sigma \backslash K)$ become pointwise
$j$-holomorphic as $n \rightarrow \infty$.  But now $j$ is smooth on
$X$, which itself is compact.  It
follows that if $|\cdbar_j(h_n)| \leq 1/n$ on $D^*$ and $h_n$ is \emph{a
  priori} smooth on the whole disc, then this bound holds
at the origin.  Each
$C_n$ defines a (much better than) rectifiable integral current.  We
assume that in fact $j$ was strictly compatible with $\omega_X$ and
write $| \cdot |$ for the associated metric.  Then
the pointwise formula

\begin{Eqn} \label{controlnorm}
|dh|^2 \ = \ h^* \omega + |\cdbar h|^2
\end{Eqn}

\noindent for any map $h: \Sigma \rightarrow X$, together with the
convergence to zero of the $|\cdbar_j (h_n)|$ on compacta of $\Sigma
\backslash K$ and the smoothness of $j$, shows that the
sequence of currents has bounded area.  We now use three results:

\begin{enumerate}
\item (Federer \cite{Federer}) The compactness theorem for bounded area
  currents guarantees a limit current $\mathcal{C}$ which will be
  \emph{area-minimizing} by (\ref{controlnorm}).  Such a current is
  necessarily a union $A \cup B$ where $B$ is the singular set and $A$
  is locally a union of finitely many submanifolds.
\item (Almgren, \cite{Almgren}) An area-minimizing integral current is
  \emph{rectifiable}; its singular set $B$ has Hausdorff codimension at
  least two.
\item (Chang, \cite{Chang})  An area-minimizing rectifiable integral
  \emph{two-dimensional} current is a classical minimal surface: it is
  the image of a smooth map $h: \Sigma \rightarrow X$ with only
  finitely many singular points.
\end{enumerate}

\noindent But now, using (\ref{controlnorm}) again, an area-minimizing
surface is necessarily $j$-holomorphic, and this completes the
argument. The caveat is that Almgren's theorem is  harder
than ``$SW=Gr$'' so even in higher dimensions there should be a better
strategy.    

\vspace{0.3cm}

\noindent $\bullet$  According to Donaldson \cite{Donaldson:pencils} the
Lefschetz pencils from pairs of  approximately holomorphic sections
become symplectic invariants when the degree $k$ of the pencils is
(arbitrarily) large.  There is a
\emph{stabilisation} procedure for Lefschetz pencils, described in
\cite{ivanmodulidivisor} and \cite{AKstabilisation}, which has the
following consequence.  Fix a degree $k$ Lefschetz pencil $\{ s_0 +
\lambda s_1 \}_{\lambda \in \pp^1}$ on $X$.  We can find a family of
degree $2k$ Lefschetz pencils, parametrised by the open complex disc
$D \subset \cc$, such that for $t \in D^*$ the pencil is a smooth
degree $2k$ Lefschetz pencil, whilst at $t = 0$ we have the degenerate
pencil $\{ s_0^2 + \lambda s_0 s_1 \}_{\lambda \in \pp^1}$.  Results 
of Auroux \cite{AurouxGokova} suggest that 
the pencils over $D^*$ do satisfy the
approximate holomorphicity constraints of \cite{Donaldson:pencils}. If we 
blow up the total space of $X \times D$ at the section of degree $4k^2
[\omega]^2$ given by the base-points of the pencils away from $0$, and
at the base curve $\{ s_0 = 0 \} \subset X$ inside $X \times \{ 0 \}$,
the resulting space carries a smooth family of Lefschetz fibrations over
$D^*$ but has singular total space.  However, further blow-ups give a
desingularisation, and one eventually obtains a smooth space with
a globally defined map to $\pp^1 \times D$.
Moreover, the fibres of this map are nodal Riemann surfaces (which may
be reducible and have many nodes each) at least away from the point
$([\lambda=0],0)$.

\vspace{0.2cm}

\noindent Now the relative Hilbert scheme can again be defined for
this larger family of curves, and at least where the singularities are
nodal and in codimension one in the base, the total space will be
smooth.  One can now hope to identify the standard surface counts for
the degree $k$ and $2k$ Lefschetz pencils with a fixed Gromov
invariant on this larger dimensional manifold (which can be
compactified or replaced with its symplectization at infinity, for
analytical purposes).  Suppose we again work with almost complex
structures making the projection to $\pp^1 \times D$ holomorphic.
Then for a generic $t \in
D$, the holomorphic curves in some class $\alpha$ will lie in a fibre
which is just the total space of a Lefschetz fibration $X'$ given by
blowing up the base-points of a degree $2k$ pencil; whilst at $t=0$,
the holomorphic curves will all lie in an irreducible component of the
relative Hilbert scheme formally obtained from a partial gluing of the
space $X_r(f)$ arising from the degree $k$ pencil and a trivial bundle
of symmetric products $\Sym^r(\{ s_0 = 0 \}) \times \pp^1$.  Hence the
moduli spaces of holomorphic sections in the total space of the larger
symmetric product fibration, for differing $J$, should recover the
spaces of sections arising from the degree $2k$ and degree $k$
pencils, setting up an equivalence between $\mathcal{I}_{(X,f_k)}$ and
$\mathcal{I}_{(X, f_{2k})}$.  Note that, in principle, the technical
difficulties here are no greater than understanding relative Hilbert
schemes for certain degenerate but integrable families of curves over
a disc.  We hope to return to this programme elsewhere.

\vspace{0.3cm}

\noindent $\bullet$  To finish this section, let us remark that the conjecture
(\ref{IequalsGr}) implies that we have a \emph{triangle} of equal
invariants:  $SW = \mathcal{I} = Gr$.  Indeed, we believe that the
equivalence of Seiberg-Witten and Gromov invariants should factor
naturally through the $\mathcal{I}$-invariants.  The link comes from the
\emph{vortex equations}, which from this point of view are just the
equations for solutions $(A, \Phi)$ to the $SW$ equations on $\Sigma
\times \rr^2$ 
which are translation invariant in the two Euclidean variables.  It is
well known that the moduli space of finite energy solutions of these
is a copy of the symmetric product of $\Sigma$.

\vspace{0.2cm}

\noindent  In \cite{Salamon}
Dietmar Salamon shows that, formally at least, the adiabatic limit of
the Seiberg-Witten equations on a fibred four-manifold, as the metric
on the fibres shrinks to zero, yields the Cauchy-Riemann equations for
holomorphic maps of the base into a universal family of solutions to
the vortex equations.  
The almost complex structure on this family arises from a canonical
connexion in
the universal vortex family.  The upshot is, formally at least, that
Seiberg-Witten solutions on $X'$ should be determined by finite energy
solutions on $X^*$, and these should come from the punctured
holomorphic sections of $\mathcal{X}^*$ that arose above.  

\vspace{0.2cm}

\noindent The vortex
equations depend on a real stability parameter $\tau$; for each fixed
value, the moduli space is a smooth symplectic manifold canonically
diffeomorphic to the geometric invariant theory quotient which is just
the usual symmetric product.  Although the complex structure is fixed, under
the adiabatic limit $\tau \rightarrow \infty$ the symplectic
structure degenerates (its cohomology class  varies linearly with
$\tau$).  Work of Hong,
Jost and Struwe \cite{Jostetal}
shows that in this limit, the solutions of the vortex equations
converge to Dirac delta solutions at the zero-sets of the Higgs field
$\Phi$.  However, the symplectic structure on the moduli space
contains a term
formally of the shape $\omega(\Phi_1, \Phi_2) = \int_{\Sigma} \langle
\Phi_1, \Phi_2 \rangle \omega_{\Sigma}$.  Thus, as Michael Hutchings
and Michael Thaddeus pointed out to the author, in the limit one
expects to obtain the degenerate symplectic structure on $\Sym^r
(\Sigma)$ given by pushing forward $\omega \times \cdots \times
\omega$ from $\Sigma \times \cdots \times \Sigma$.  This is analogous
to degenerating almost complex structures on $X_r(f)$ to ones coming
from the fibre product.


\section{Serre Duality for symplectic surfaces}

In this section we shall prove Theorem
(\ref{Serreduality}).  As we observed in the Introduction, the
geometric input for the first part arises from Serre duality on the
fibres of a Lefschetz 
fibration.  In order to compute Gromov invariants, however, we need
almost complex structures which behave well when we pass to
one-parameter families of curves.  The key
technical ingredient allowing this is a result from Brill-Noether
theory, due to Eisenbud and Harris, which we introduce at once.  


\subsection{Brill-Noether theory}

Let $C$ be a Riemann surface of genus $g$.
Let $W^s_r(C)$ denote the locus of linear systems
$g^s_r$ of degree $r$ and dimension $s$ on $C$, viewed as a subscheme
of the Picard variety $\Pic^r(C)$.  This has virtual dimension
given by the \emph{Brill-Noether number}

$$\rho  \ = \ g-(s+1)(g-r+s).$$

\noindent The famous ``existence theorem'' \cite{ACGH} asserts that if
$\rho \geq 0$ then $W^s_r(C)$ is not empty, for \emph{any} $C \in
M_g$.  On the other hand, the ``dimension theorem'' (\emph{op. cit.})
asserts that if $\rho < 0$ then $W^s_r(C)$ is empty for a
\emph{general} curve of genus $g$; in other words, it is non-empty
only on a subvariety of positive codimension in $M_g$.  One could hope to
sharpen this statement, at least in good cases, to estimate the
codimension of this subvariety in terms of the deficiency $-|\rho|$ of the
(negative) Brill-Noether number.  In general the naive estimates fail (as shown
in many cases in \cite{ACGH}), but eliminating codimension one
components does turn out to be possible \cite{EH}:

\begin{Thm}[Eisenbud-Harris] \label{EHresult}
Suppose $\rho < -1$.  Then the locus $\mathcal{W}^s_r \subset M_g$
comprising curves $C$ for which $W^s_r (C) \neq \phi$ has codimension
greater than one in the moduli space.
\end{Thm}

\noindent In \cite{SKDIS} a
central ingredient in the computation of the Gromov invariant was the
simple form of the Abel-Jacobi map $X_{2g-2}(f) \rightarrow
P_{2g-2}(f)$, with a unique fibre of excess dimension.  This holds for
any choice of fibrewise metrics on the Lefschetz pencil.  Such a simple
description cannot exist for all Abel-Jacobi maps $X_r(f) \rightarrow
P_r(f)$; if $r=g-1$ then the geometry of the theta-divisor is well
known to be subtly dependent on the curve in moduli space.  However, using the
above proposition and the adjunction formula, we can obtain a pretty
description for the cases of interest to us.  Fix a symplectic
four-manifold $X$ and let $\alpha \in H_2(X;\zz)$.  Under the map
$\iota: H_2(X; \zz) \rightarrow H_2(X'; \zz)$ we map $\alpha$ to a
class $\alpha + \sum_{i=1}^N E_i$ whose intersection with the fibre of the
Lefschetz fibration $f: X' \rightarrow \sss^2$ is at least as large as
the number $N$ of exceptional curves.  But by adjunction on $X$, if the
fibres of the pencil represent a class $W$ in homology,

$$2g(\mathrm{Fibre}) - 2 \ = \ K_{X} \cdot W + W^2,$$

\noindent where of course $N = W^2$.  As $W = k[\omega]$ then by
increasing the degree of the 
pencil, the second term on the right grows quadratically with $k$ and
the first term on the right only linearly; hence the ratio $(2g-2) / N$
tends to $1$. In particular, for any $\alpha$, the intersection number
$r$ of $\iota(\alpha)$ and the fibre of the Lefschetz fibration
is such that $(2g-2)/r \rightarrow 1$ as the degree increases.

\begin{Prop} \label{nicefibres}
For sufficiently large $k$ and $r(k)$, there is an embedding of the
residual fibration $T \cong X_{2g-2-r}(f) \hookrightarrow P_r(f)$.  The
natural map  $u:X_r(f) \rightarrow P_r(f)$ has (projective space) fibres
of dimension $r-g$ away from $T$ and of dimension $r-g+1$ over $T$.
\end{Prop}

\begin{Pf}
By the above remarks, we can certainly assume that $r>g-1$. Fix a
smooth curve $C$.
The fibre of $u:\Sym^r (C) \rightarrow \Pic^r (C)$ has dimension $r-g+1$
for line bundles with no higher cohomology.  Conversely, the line
bundle $L$ has higher cohomology if and only if $K-L$ admits sections;
such a section defines an effective divisor of degree $2g-2-r$, and
hence the space $\Sym^{2g-2-r}(C)$ parametrises all of the residual
linear systems at which points the map $u$ may have vanishing differential.
Now for $2g-2 \geq r \gg g$, we see $\rho$ is
negative: setting $d=2g-2-r$ to be the degree of the
residual divisors, we have 

$$\rho(W^1_d(C)) = g-2(g-1+d) = 2-g-2d \leq 2-g$$

\noindent which will be substantially smaller than $-1$ for high
degree pencils.  The virtual dimensions of the $W^s_d(C)$ with $s>1$
are still smaller; since $d < g-1$ we expect the line bundle $K-L$ to
have no sections, and it to be increasingly unlikely that linear
systems of dimension $s$ exist as $s \mapsto s+1$.  Hence for a
generic $C$, meaning a point $C \in M_g$ in the complement of a
subvariety of codimension at least $2$, the Abel-Jacobi map
$\Sym^r(C) \rightarrow \Pic^r (C)$ has excess fibres of dimension at
most one greater than the generic dimension, and moreover the
locus of points in $\Pic^r(C)$ where this happens is exactly a copy of
the dual symmetric product.  This embeds in $\Pic^r(C)$ under the
obvious map $\sum z_i \mapsto (K_C-\mathcal{O}(\sum
z_i))$.

\vspace{0.2cm}

\noindent It follows, by the
Eisenbud-Harris result, that for a generic choice of fibrewise
metrics the sphere $\phi_f (\sss^2) \subset \mgbar$ will be disjoint
from the locus $\mathcal{W}^s_r$ for every $s$, at least over the
smooth locus.  We now need to argue for the behaviour over the
singular fibres.  Recall from \cite{Altman-Kleiman} that there is always
a natural morphism $\Hilb^{[d]}(\Sigma_0) \rightarrow \Pic_d
(\Sigma_0)$, for every $d$, with projective space fibres;  indeed, a
generic point of $\Pic_d$ gives a locally free sheaf on the
normalisation and the projective space is
just the linear system of divisors of this line bundle.  On the moduli
space $M_{g-1}$ there is a family of complex surfaces given by the
fibre product of the universal curve with itself, and a natural
(gluing) morphism from this onto the divisor in $M_g$ of irreducible
curves with one node.  The Eisenbud-Harris result, applied for genus
$g-1$ curves, now shows that the embedding we have constructed away
from the singular fibres will extend to embed $T \hookrightarrow
P_r(f)$.  The rest of the result is clear; for every
degree $r$ torsion-free sheaf $L$ on each (not necessarily smooth) curve $C$
in the Lefschetz 
fibration, the dimension of the linear system $K-L$ is at most zero,
meaning that the effective divisors of degree $2g-2-r$ never move in
non-trivial systems. 
\end{Pf}
\begin{Example} As in the Appendix to
  \cite{SKDIS} the situation is most clear for the case of divisors
  of degree two.  
Let $\Sigma_0$ be a nodal curve of high genus; the Hilbert scheme
$\Hilb^{[2]}(\Sigma_0)$ is given by blowing up the second symmetric
product at the point $(\mathrm{Node}, \mathrm{Node})$.  $\Pic_2
(\Sigma_0)$ is given by gluing two sections of a $\pp^1$-bundle over
$\Pic_2 (\tilde{\Sigma}_0)$ over a translation in the base.  The
embedding $\Hilb \rightarrow \Pic$ defines a map from $\Hilb^{[2]}$ to
the normalisation, which maps the exceptional sphere to a fibre.
Inside $\Pic_2$ the sphere remains embedded, and meets the singular
locus at two points, which are just the two points of the exceptional
curve singular inside $\Hilb^{[2]}$.  
\end{Example}


\subsection{Adapted complex structures}

We now have a very simple picture in which to compute the
standard surface counts of homology classes.  Suppose we have a fibre
bundle $\pi: Z \rightarrow B$ with
almost complex fibres $F_b$ and almost complex base $B$.  Any choice
of connexion (horizontal splitting) for $\pi$ defines an almost
complex structure on $Z$; we use the given $J_F$ on the tangent spaces
vertically and use the connexion to lift $j_B$ to the horizontal
planes.  Now connexions may always be extended
from closed subsets - their obstruction theory is trivial.  (As usual:
given \emph{some} connexion on the fibration, a
particular connexion on a subfibration is given in each point by the
graph of a linear map from the horizontal to the vertical space at
that point.  Since the space of such linear homomorphisms is
contractible, the given section over the closed subspace may be
extended.)  It follows
that we may find almost complex structures on the total space making
any subfibration by almost complex subspaces of the fibres an almost
complex submanifold of the total space.

\begin{Defn}
In the situation of the previous proposition, there are almost complex
structures on $P_r(f)$ and $X_r(f)$ extending the obvious integrable
structures on the fibres and for which
\begin{enumerate}
\item the submanifold $X_{2g-2-r}(f) \rightarrow P_r(f)$ is an almost
complex subspace;
\item the Abel-Jacobi map $u:X_r(f) \rightarrow P_r(f)$ is holomorphic;
\item the almost complex structure on the restriction $u^{-1}
  (X_{2g-2-r}(f))$ is induced by a linear connexion on the total space
  of a vector bundle over $X_{2g-2-r}(f)$.
\end{enumerate}
We will call such almost complex structures $J \in \mathcal{J}$ on $X_r(f)$
\emph{compatible with duality}.  If $r=2g-2$, as in \cite{SKDIS}, we
will say $J$ is \emph{compatible with the zero-sections}.
\end{Defn}

\noindent This is an immediate consequence of the discussion.  Choose
an almost complex structure on $X_{2g-2-r}(f)$ and then a
connexion on a vector bundle over this space with fibre at a divisor
$D$ the sections $H^0(K-\mathcal{O}(D))$.  That these complex vector spaces,
by assumption now of constant rank, fit together to yield a vector
bundle is essentially standard elliptic regularity (cf. the next
proposition).  More easily,
the vector bundle is just the pullback by $\phi_f: \sss^2 \rightarrow
\mgbar$ of a universal vector bundle over the moduli space.  This is
the push-forward of the bundle over $\mathcal{S}^{2g-2-r}(\mathcal{C})
= \Hilb^{[2g-2-r]}(\mathcal{C}_g /
\mgbar)$ with fibre at a subscheme $Z \subset C$ the line bundle $K_C
- \mathcal{O}(Z)$.  The ambiguity arising from the fact that the
relative dualising sheaf of a fibration is not the same as the
restriction of the canonical bundle of the total space to each fibre
is eliminated since it is only the \emph{projectivisation} of the
vector bundle which embeds into $X_r(f)$.  Now given a complex
structure on the total space of the projective bundle, viewed as a
closed subset inside $X_r(f)$, we can extend arbitrarily by choosing
an extension of the connexion to the whole of $u: X_r(f) \rightarrow P_r(f)$.

\vspace{0.2cm}

\noindent Since these vector bundles of sections on the fibres of a
Lefschetz fibration appear as holomorphic subspaces of the $X_r(f)$,
it will be helpful to understand their spaces of holomorphic sections.
Here is a general result in this line, formulated
in the framework and notation of Lemma (\ref{Cyclesinject}).
Let $\alpha \in H_2(X;\zz)$ and bear in mind (\ref{samedim}).

\begin{Prop} \label{indextheorem} Suppose, for each $j \in \{ 0,1 \}$,
  the value  $h^j
  (L_{\alpha}|_{\Sigma_t})$ is constant for $t \in \sss^2$.  Then the
  vector spaces $H^0 (L_{\alpha}|_{\Sigma_t})$ and $H^1
(L_{\alpha}|_{\Sigma_t})$ define vector bundles over $\sss^2$. The 
element $V_{\alpha} = \{ H^0 - H^1 \}_{t \in \pp^1}$ of the $K$-theory
  of $\cc \pp^1$ satisfies

$$\ind(V_{\alpha}) \ = \ [\alpha^2 - K_X \cdot \alpha]/2 + (b_+ +
1-b_1)/2.$$

\noindent Here the index $\ind(\cdot)$ is defined as the sum of the
rank and first Chern class of the virtual bundle.
\end{Prop}

\begin{Pf}
The existence of the vector bundles is an application of
standard results in index theory for families of elliptic operators.
Away from the critical fibres we have a smoothly varying family of
$\cdbar$-operators.  Over small discs around each critical value of the
Lefschetz fibration, we can assume both $X$ and the families of
operators vary holomorphically, and the index virtual bundle is just the
push-forward $f_! (L_{\alpha})$.  There is a canonical identification
of these two objects over annuli around each critical value, coming
from the canonical holomorphic structure on the index bundle for a
family of Hermitian operators on the fibres of a holomorphic submersion,
as explained in Freed's survey \cite{Freed}.

\vspace{0.2cm}

\noindent The first Chern class of this index bundle could be computed
using either the Atiyah-Singer theorem, were the family smooth (no
critical fibres), or the
Grothendieck-Riemann-Roch theorem were it globally holomorphic (even with
critical fibres).  Although neither apply here, we can get around the
deficiency.  One approach (which we shall not develop) uses excision,
and removes a neighbourhood of the critical fibres.  Alternatively, our
family of $\cdbar$-operators is classified by a smooth map $\psi:
\sss^2 \rightarrow \Pic$ to a universal family over the universal
curve.  The homotopy invariance property for indices of families of
operators asserts that the index is determined by the homotopy class
of this map, which is encoded by the homology class $\alpha$.  The
index bundle $V_{\alpha}$ is just the pullback $\psi^*
(\pi_! \mathcal{F})$ of the push-forward of a holomorphic sheaf
$\mathcal{F} \rightarrow \Pic$.  We can now argue precisely as in the proof of
(\ref{virdim}) to obtain the formula.
\end{Pf}

\noindent Note that in the proof of (\ref{virdim}) and above, the
universal perspective is used to reduce the computation to the
specific case of the first Chern class of the relative dualising
sheaf, which is computed separately in \cite{ivanhodge}.  It would be
interesting to develop a direct argument for giving the indices of
smooth families of operators with locally holomorphic singularities.

\begin{Rmk} \label{exactseqs}
Let us stress again the two key geometric features of both the index
formulae we have given: the same geometry, in a mildly more
complicated scenario,  will play a role in the next
subsection.  For complex structures compatible with the zero-sections, the
geometry of the situation is reflected in maps:

$$\pp(V) \hookrightarrow X_r(f) \stackrel{F}{\longrightarrow} P_r(f).$$

\noindent These induce:

$$0 \rightarrow T\pp(V) \rightarrow T^{vt} X_r(f) \rightarrow \im (dF)
\rightarrow 0;$$
$$0 \rightarrow \im (dF) \rightarrow T^{vt} P_r(f) \rightarrow \cok
(dF) \rightarrow 0.$$

\noindent  If we take the two long exact sequences in cohomology, we
obtain a sequence:

$$0 \rightarrow H^0 (W) \rightarrow \mathrm{Obs}_u \rightarrow
\cc^{(b_+-1)/2} \rightarrow H^1 (W) \rightarrow 0.$$

\noindent Here we have used two identifications:

\begin{enumerate}
\item At any point $\sum p_i \in \Sym^r(\Sigma)$ the cokernel
of the Abel-Jacobi map is canonically isomorphic to $H^1 (\Sigma,
\mathcal{O}(\sum p_i))$ and this globalises to the identity $\cok (dF)
\cong W$, where $W$ denotes the appropriate push-forward sheaf.  Now
the first and last terms in the sequence are
just $H^i (\cok (dF)_u)$ for $i \in \{0,1\}$.
\item The vertical tangent bundle to $P_r(f)$ near any section, by the
  remarks above, is isomorphic to the bundle with fibre $H^1 (\Sigma,
  \mathcal{O})$.  Relative duality - a parametrised form of Serre
  duality which holds for any flat family of curves, and hence for the
  universal curve from which our bundles are pulled back
  \cite{HartshorneRD} - asserts that

$$T^{vt} P_r(f)|_{u(\sss^2)} \, \cong \, (f_* \omega_{X'/\pp^1})^*. $$
\end{enumerate}
\end{Rmk}

\noindent These comments on vertical tangent bundles implicitly assume
that all our sections are smooth, in the sense that there are no
bubbles.  We have proved that fact for regular almost complex structures,
and shown that structures compatible with the diagonal strata can be
regular, but we are now working with almost complex structures which
fall into neither of these classes.  Hence we must provide a fresh
argument.

\begin{Prop}
Let $J \in \mathcal{J}$ on $X_r(f)$ be generic amongst structures which are
compatible with duality.  Then
for any holomorphic section $w$ of $X_{2g-2-r}(f)$, the holomorphic
sections of the 
projective subbundle of $X_r(f)$ lying over $w$  arise from constant
sections of a vector bundle and contain no bubbles.
\end{Prop}

\begin{Pf}
Fix the section $w$ of $X_{2g-2-r}(f)$.  We will work in the preimage of
this fixed section.
The holomorphic sections of a projective bundle $\pp(V)$ coming from
sections of $V$ have bubbles iff the sections of $V$ vanish at certain
points.  For us, once we have fixed a homology class $\alpha \in
H_2(X;\zz)$, 
the relevant vector bundle $V^0$ is given by taking the holomorphic
sections of the bundles $L_t = L_{\alpha}|_{\Sigma_t}$ on the fibres
of the Lefschetz 
fibration.  We have already given the index of the $\cdbar$-operator on the
element $V^0 - V^1$ in K-theory, where $V^1$ is the bundle of first
cohomology groups on the fibres.  For high degree pencils the
rank of $V^0$ is very large, growing with $k$, whilst $J$ being
compatible with duality precisely asserts that the bundle $V^1$ has
rank one.  Hence it is enough to show that the first Chern class of
$V^1$ is small, and then the constant value of $\mathrm{Index}(V^0 -
V^1)$ will force the first Chern class of $V^0$ to be negative for
high degree pencils.  In this case, by stability, the bundle will
generically be of the form

$$V^0 \ = \ \mathcal{O} \oplus \cdots \oplus \mathcal{O} \oplus
\mathcal{O}(-1) \oplus \cdots \oplus \mathcal{O}(-1);$$

\noindent the only holomorphic sections of such a bundle over $\pp^1$
are constant sections, and these have no isolated zeroes and give no bubbles.

\vspace{0.2cm}

\noindent In fact, we claim that the first Chern class of the line
bundle $V^1$ is necessarily zero, in other words the line bundle is
topologically trivial.  (In the situation of \cite{SKDIS} this amounts
to saying that the
push-forward sheaf $R^1 f_* \mathcal{O}$, down the fibres of a
Lefschetz fibration, is trivial.)  The result is an application
of relative duality.  Fibrewise, $H^1 (\Sigma_t;
L_{\alpha}|_{\Sigma_t})$ is dual to $W_t = H^0 (K_{\Sigma_t} -
L_{\alpha}|_{\Sigma_t})$; hence it is enough to understand the Chern
class of the bundle with this as fibre.  But a section $u$ of $W
\rightarrow \pp^1$ 
gives, at each value $t$, a section $u_t$ of the residual bundle and
hence a divisor in $\Sym^{2g-2-r}(\Sigma_t)$ - which is necessarily
the unique point $w(t)$ in the linear system, under our assumption that
$\Sym^{2g-2-r} \hookrightarrow \Pic_r$.
There is a subtlety here;  the section of the projective bundle is
only a divisor on each fibre up to scale, and so one might think this
scalar indeterminacy precisely hid the first Chern class.  In fact,
\emph{if the section $u$ has a zero, then its zero-set will contain
  certain fibres}, but the collection of points of the symmetric
products defined by the original section $w = \pp(u)$ does not.  

\vspace{0.2cm}

\noindent More invariantly, argue as follows.  In general, we can't
assume that an almost complex structure $J$
compatible with duality also respects the stratification coming from
exceptional sections (i.e. makes the subfibrations $X_{r-1}^E(f)$ all
holomorphic).  However, if we distinguish one fixed $E$, we can also
assume $X_{r-1}^E(f)$ is an almost complex subvariety, since one
condition on a vector space is always linearly independent and cuts
out a hyperplane.  By the usual arguments that enable us to push
cycles from $X'$ down to $X$, the original section $w$ of
$X_{2g-2-r}(f)$ must in fact be \emph{disjoint} from $X_{r-1}^E(f)$.
It follows that, by scaling to make the evaluation of the section $u$
at the unique point of $E$ in any given fibre equal to $1$, we can
explicitly smoothly trivialise $W$.  This completes the argument.
\end{Pf}

\noindent A consequence of the above result, and the preceding
remarks, is the following: the arguments are similar so we leave them
to the reader.

\begin{Cor} 
The normal bundle of the embedding of $X_{2g-2-r}(f)$
inside $P_r(f)$ has negative first Chern class. 
\end{Cor}

\noindent With these two facts in hand, we can assemble all the pieces
for our version of Taubes' duality theorem.


\subsection{Proof of the duality}

The proof will be completed as follows.  We shall fix an
almost complex structure on $X_r(f)$ which is compatible with duality,
and explicitly determine the moduli space of smooth holomorphic sections.
Despite the non-genericity, this moduli space will already be compact,
and hence we will be able to compute the obstruction bundle (by
using elementary arguments to understand its fibre at any given smooth
section). 
For the moment, let's work with an easier linear constraint:

\begin{Thm} \label{dualityproof}
Let $X$ be a symplectic four-manifold with $b_+ (X) > 1+b_1 (X)$.
For any class $\alpha \in H_2 (X; \zz)$ the standard surface counts
for $\alpha$ and $\kappa - \alpha$ co-incide up to sign.
\end{Thm}

\begin{Pf}
By the blow-up formula (\ref{blowup}) given before, 
it is enough to compare the counts of holomorphic sections giving curves
in classes $\iota(\alpha)$ and $\kappa_{X'} - \iota(\alpha) = i(\kappa_X
- \alpha)$.   

\vspace{0.2cm}

\noindent For $\iota(\alpha)$, the algebraic intersection number with the
fibre of $f: X' \rightarrow \sss^2$ is $r \gg g-1$ and we can assume
that the fibration $X_r(f)$ maps to $P_r(f)$ with the standard
topological format
described above.  Fix an almost complex structure $J \in \mathcal{J}$ which
is generic on $X_{2g-2-r}(f)$ and which is generic amongst those
compatible with duality in
the sense of the preceding Proposition.  Since the Abel-Jacobi map is
assumed holomorphic, every ${J}$-holomorphic section of $F$
maps to a holomorphic section of $P_r(f)$.  We claim that, for
our particular ${J}$, this must
lie inside $X_{2g-2-r}(f) \subset P_r(f)$.  To see this, recall that the
index of the $\cdbar$-operator on the vector bundle $T^{vt}P_r(f)
\rightarrow \sss^2$ is negative, more precisely $1+b_1-b_+$.  This
follows since for any fixed section $\phi$ of
$P_r(f)$ there is a diffeomorphism of pairs 

$$(P_r(f), \phi) \rightarrow (P_{2g-2}(f), \phi_K)$$

\noindent where $\phi_K$ is the section coming from the canonical
bundles on the fibres.  This identifies the vertical tangent bundle
near $\phi$ with that near $\phi_K$, and the first Chern class of this
bundle determines the index:  but then the index for the section
$\phi_K$ is determined by the signature formula of \cite{ivanhodge},
as explained earlier in the paper.  (This is the $\alpha = 0$ instance
of \ref{indextheorem}.)

\vspace{0.2cm}

\noindent It follows that for any holomorphic section $\phi$ of
$P_r(f)$ passing through some (previously) fixed open subset $U
\subset P_r(f) \backslash X_{2g-2-r}(f)$, the moduli space of
holomorphic sections in the class
$[\phi]$ is regular near $\phi$ and hence of the correct (virtual)
dimension, which is negative.  We observed above that the normal bundle $N$
to the embedding of
$X_{2g-2-r}(f)$ inside $P_r(f)$ has negative Chern class.  Therefore if $J$ is
generic on a
tubular neighbourhood of the image of the embedding, and vertically
generic on its preimage under the Abel-Jacobi map, 
we can deduce that in fact $N$ has no
sections either.  Hence sections of $X_{2g-2-r}(f)$ do not deform
infinitesimally in the total space unless they remain inside this
subfibration.  This means that all holomorphic sections for the
almost complex structure ${J}$ have image inside the embedded copy of
$X_{2g-2-r}(f)$.   By the
definition of the embedding 

$$X_{2g-2-r}(f) \ \rightarrow \ P_r(f); \qquad \sum z_i \ \mapsto \ K
- \mathcal{O}(\sum z_i)$$

\noindent the section $u \circ \phi: \sss^2 \rightarrow X_{2g-2-r}(f)$
gives rise to a cycle in the homology class $i(\kappa - \alpha)$ if
$\phi$ gives a cycle in the class $\iota(\alpha)$.  Conversely, for
any holomorphic section of $X_{2g-2-r}(f)$ in the former class, the
restriction of $u$ to the preimage of this section is just the
projection map of the projectivisation of a vector bundle $V$, by the
result (\ref{nicefibres}).  Hence for some $a=\ind_V(\cdbar)-1$, we have
exhibited the moduli space of ${J}$-holomorphic sections of
$X_r(f)$ in the class $\phi^{-1}(\iota(\alpha))$ as a $\pp^a$-bundle over the
moduli space of ${J}$-holomorphic sections of $X_{2g-2-r}(f)$ in the
class $\phi^{-1}(i(\kappa - \alpha))$.   It remains to compare the
obstruction bundles on the two moduli spaces.  

\vspace{0.2cm}

\noindent First, since we have a generic almost complex structure
downstairs, we can choose $m = [\alpha^2 - K\cdot\alpha]/2$ points on
$X$ and cut down the moduli space of sections of $X_{2g-2-r}(f)$ to be
zero-dimensional, and hence a signed set of points.  Their total
number is exactly the standard surface count for
$\kappa-\alpha$.  Since the rank of the linear system over a
point of $P_r(f)$ is large relative to $m$, each condition cuts down
each linear system by a hyperplane (the conditions are independent).
Thus the new geometric situation is again that of a moduli space which
is a projective bundle but now over finitely many points.  The
obstruction theories behave compatibly,
so we reduce to studying the $m=0$ problem for a smaller family of
symmetric products $\Sym^{r-m}$; so assume $m=0$ from the start.  
It will be enough to show that the obstruction bundle on
each $\pp^a$ is just the quotient bundle (that is, the cokernel of the
embedding of the tautological line bundle into the trivial bundle of
rank $a+1$).  Then each $\pp^a$ will
contribute $\pm 1$ to the moduli space of sections in the class
$\iota(\alpha)$, since the quotient bundle has Euler class $\pm 1$
depending only on the dimension $a$.

\vspace{0.2cm}

\noindent Write $Z$ for the total preimage of $X_{2g-2-r}(f)$ inside
$X_r(f)$ under the Abel-Jacobi map.  Recall that the
sections we are interested in stay away from all the critical loci of
the projection maps, and so we can consider the cokernel $W$ of the
natural projection $\overline{\tau}$:

\begin{Eqn} \label{cokernel}
0 \rightarrow \ T(X_r(f))\big{/}TZ \ \
\stackrel{\overline{\tau}}{\longrightarrow} \ \ 
T(P_r(f))\big{/}T(X_{2g-2-r}(f)) \ \longrightarrow W \ \rightarrow 0.
\end{Eqn}

\noindent We also have, for a given section $s \in \pp^a$ in one
component of our moduli space, a sequence

$$0 \ \rightarrow \ \nu_{s/Z} \ \rightarrow \ \nu_{s/X_r(f)} \
\rightarrow \ \nu_{Z/X_r(f)}|_s \ \rightarrow 0$$

\noindent of normal bundles.  Here we identify $s$ with its image, and
the sequence is induced by the inclusions $s \subset Z
\subset X_r(f)$.  The long exact sequence in cohomology for this
second sequence shows

\begin{Eqn} \label{obstructionsequence}
0 \ = \ H^1 (\nu_{s/Z}) \ \rightarrow \ H^1(\nu_{s/X_r(f)} = Obs(s)
\ \rightarrow \ H^1 (\nu_{Z/X_r(f)}|_s) \  \rightarrow \ 0.
\end{Eqn}

\noindent The first identity holds since the almost complex structure
${J}$ is generic over the \emph{whole} of $Z$; the
identification in the middle is just the definition of the obstruction
space at $s$.  Now we use the sequence (\ref{cokernel}):

\begin{Eqn} \label{laststep}
0 \rightarrow H^0 (W) \rightarrow Obs(s) \rightarrow H^1 \big
( T(P_r(f))\big{/}T(X_{2g-2-r}(f)) \big)|_s \rightarrow H^1(W)
\rightarrow 0.
\end{Eqn}

\noindent The penultimate term depends only on the image of $s$ in
$P_r(f)$ so is constant as we vary over the fixed projective space;
i.e. it varies to give a trivial bundle over $\pp^a$.  We must
identify the cokernel bundle $W$.  The fibre
of $W$ at a divisor $D = \sum p_i$ is exactly (canonically)
$H^1 (\Sigma;D) = H^0 (\Sigma; K-D)^*$ \cite{ACGH}.  On the other hand, $H^0
(K-D)$ is \emph{generated by the unique effective divisor} which is
the point $u(D) \in X_{2g-2-r}(f)$, where $u$ is the reciprocity map.
It follows, as in \cite{SKDIS}, that we can identify $W \cong
K_{\pp^1}$, where the twist comes from the relative duality
isomorphism, and then the fact that $u(D)$ generates $H^0 (K-D)$
identifies the sequence with:

$$0 \rightarrow \mathcal{O}(-1) \rightarrow \cc^{a+1} \rightarrow
\mathrm{Obs} \rightarrow 0.$$

\noindent The first term is the tautological bundle, and hence the
(quotient) obstruction bundle is the quotient bundle on projective space
as claimed.

\vspace{0.2cm}

\noindent The above applies to each component $\pp^a$ of the moduli
space of sections representing $\iota(\alpha)$.  But the Euler class
of the quotient bundle is $\pm 1$ where the sign is determined by
$a$.  This is fixed once and for all by $\alpha$ and $X$, so we indeed
find that the standard surface counts for $\alpha$ and $\kappa-\alpha$
can only differ by a single overall sign.
\end{Pf}

\noindent As an immediate application, here is the ``simple type''
result that we alluded to earlier.

\begin{Prop}
Let $X$ be a symplectic four-manifold with $b_+ (X) > 1+b_1(X)$.  Then
the standard surface counts vanish for any class $\alpha$ with
$\alpha^2 \neq \kappa \cdot \alpha$.
\end{Prop}

\begin{Pf}  Fix a
compatible almost complex structure $J$ on $X$. Suppose the standard
surface count for a class $\alpha \in H_2 (X; \zz)$ is non-zero.
Then, by the symmetry proven above and the results of
the previous section,
both $\alpha$ and $\kappa-\alpha$ can be represented by embedded (or
smoothly multiply covered) $J$-holomorphic 
curves in $X$.  It is well known that $J$-holomorphic curves have
locally positive intersections in a four-manifold;  hence $\alpha$ and
$\kappa-\alpha$ must have non-negative intersection unless they have
common components.  Suppose $C$ is such a component.  Then $\alpha
\cdot C = C^2 = (\kappa-\alpha)\cdot C$, which means that $\kappa
\cdot C = 2C^2$.  On the other hand, by adjunction, $\kappa \cdot C +
C^2 = 2g(C)-2$ and hence $C^2 \geq 0$.  It follows that in all cases
$\alpha \cdot (\kappa - \alpha) \geq 0$.
Contrastingly, the (real) dimension of the moduli space of holomorphic
sections of $X_r(f)$ in the class corresponding to $\alpha$ is exactly
$\alpha^2 - \kappa\cdot\alpha$, also proven above; this must be
non-negative if the
invariant is non-zero.  Hence

\begin{Eqn} \label{firstset}
\alpha^2 \leq \kappa \cdot \alpha; \qquad \alpha^2 \geq \kappa \cdot
\alpha.
\end{Eqn}

\noindent This gives $\alpha^2 = \kappa \cdot \alpha$, as required.
\end{Pf}

\noindent It's now easy to check the following, which gives one
familiar ``basic class'' type obstruction to symplectic manifolds
being K\"ahler.  

\begin{Cor}
Suppose $b_+ (X) > 1+b_1 (X)$.
If $\mathcal{I}_{(X,f)}(\alpha) \neq 0$ then $0 \leq \alpha \cdot
[\omega] \leq K_X \cdot [\omega]$, with equality iff $\alpha \in \{0,
K_X \}$.  If in fact $(X;\omega)$ is a minimal K\"ahler surface of general
type, then the standard surface counts are non-zero iff $\alpha \in 
\{ 0, K_X \}$.
\end{Cor}

\begin{Pf}
The first half is immediate from the above.  The result on the
invariants for general type surfaces follows from the fact that $K_X$
is in the closure of the ample cone.
Indeed, using results on pluricanonical linear systems \cite{BPV},
it's easy to see
that in fact $K_X$ contains symplectic forms deformation equivalent to
the given $\omega$ (just take an exact perturbation near any
contracted $(-2)$-spheres).  So then if $\mathcal{I}_{(X,f)}(\alpha)
\neq 0$ then 

$$\alpha \cdot [\omega] \leq K_X \cdot [\omega] \ \Rightarrow \ \alpha
\cdot K_X \leq  K_X^2.$$

\noindent Now the Hodge index theorem and Cauchy-Schwartz assert that
$(\alpha^2) \cdot (K_X^2) \leq (\alpha \cdot K_X)^2$, which forces
either $\alpha = 0$ or an equality $\alpha^2 = K_X^2$.  The result is
a simple consequence. 
\end{Pf}


\subsection{More symplectic surfaces}

In this subsection, we prove the second part of Theorem
(\ref{Serreduality}).  When $b_+ = 1$ the virtual dimension for
sections of the Picard fibration is \emph{never} negative, and so the
arguments used above to control the geometry of the moduli spaces must
be refined somewhat.

\begin{Lem} \label{torusofsections}
Let $(X,f)$ be any symplectic Lefschetz pencil on a manifold with
$b_+ = 1$.  Form the Picard bundle
$P_r(f) \rightarrow \sss^2$ and fix a section $s$ of $P_r(f)$.  Then
there is an almost complex structure $j$
on $P_r(f)$ for which the moduli space of sections in the homotopy
class defined by $s$ is exactly a complex torus $\rot^{b_1(X)/2}$.
\end{Lem}

\begin{Pf}
Note first the statement is coherent;  if $b_+=1$ then $b_1$ is even.
Decompose the base $\sss^2$ into a union $D$ of discs around the
critical values of $f$ over which the fibration is holomorphic, annuli
$A$ surrounding $D$ and the complement $R$ of slightly smaller annuli.
Over $R$ the degree zero Picard fibration (that is, the Jacobian
fibration) carries a canonical flat connexion
defined by the integral lattices in each Picard torus fibre, and this
defines an almost complex structure on the total space which is linear
in the sense that it lifts to a linear connexion on the complex vector bundle
$V|_R \rightarrow R$ where $V \rightarrow \pp^1$ has fibre $H^1
(\Sigma_t, \mathcal{O})$ at $t \in \pp^1$.  Over the discs $D$ the
total space of the Picard bundle is exactly isomorphic to the sheaf
quotient $R^1 f_* \mathcal{O} / R^1 f_* \zz$ and hence the induced
integrable complex structure on $P_r(f)|_D$ lifts to the given
linear complex structure on the same vector bundle $V$.  Now over the
interpolating annuli $A$, we can interpolate the two connexions on $V$
equivariantly with respect to the action of the lattice $\zz^{2g}$
(fix a homotopy on a fundamental domain and then extend).
This defines a complex structure on $P_0(f)$.

\vspace{0.2cm}

\noindent Now be given some $P_r(f)$ and a fixed section $s$.  We use
the section to define a diffeomorphism of $P_r(f)$ and $P_0(f)$ which
takes $s$ to the zero-section, and equip $P_r(f)$ with the pullback of
the above complex structure.  Clearly the given section is now
holomorphic.  Given any other holomorphic section $s'$ of
$P_r(f)$ which is homotopic to $s$, the differences $y(t) = s(t)^{-1} \otimes
s'(t)$ define a section of the bundle $P_0(f)$ which is also
holomorphic.  Since we have worked with integrable structures near the
singular fibres, and the group action is holomorphic in the algebraic
setting, this is globally well-defined.  By the construction of the
complex structure from a
connexion, this section $y$ lifts (non-uniquely!) to a section of the
vector bundle $V
\rightarrow \pp^1$, which for a generic interpolating choice of
connexions over the annuli $A$ carries its most stable complex
structure.  But the index of the $\cdbar$-operator on $V$ is just
$b_1(X)/2$ and the rank of the bundle is $g \gg 0$. Hence the most
stable complex structure is of the shape 

$$V \cong \mathcal{O} \oplus
\cdots \mathcal{O} \oplus \mathcal{O}(-1) \oplus \cdots \oplus
\mathcal{O}(-1)$$

\noindent and the moduli space of sections of $V$ in the required
homotopy class is just $\cc^{b_1(X)/2}$.  To compute the moduli space
of sections of $P_r(f)$ in the projected homotopy class, we must
divide out the non-uniqueness of the lift, and this is just the action
of the sublattice that preserves the holomorphically trivial subbundle of
$V$.  This yields a torus of the claimed dimension.
\end{Pf}

\noindent Let us put this observation in context.  Suppose for a
moment that we are looking for symplectic surfaces in a class $\alpha$
which arise from a section of the bundle of symmetric products
$X_r(f)$ with $r>2g-2$.   This is always the case, by adjunction, for a
class of the form $\iota(\alpha)$ when $K_X \cdot \omega < 0$.  In
this case the geometry
simplifies considerably;  the total space of $X_r(f)$ is a projective
bundle over $P_r(f)$.  The map $F: X_r(f) \rightarrow P_r(f)$ is a
submersion away from the critical fibres, and the exact sequences of
(\ref{exactseqs}) and (\ref{dualityproof}) disappear.  We find that
\emph{the obstruction bundle is topologically trivial}.  This has two
consequences:

\begin{itemize}
\item If $b_+ > 1$ then the obstruction bundle has trivial Euler
  class, and hence the standard surface counts cannot be non-trivial
  for any class $\alpha$ with $\alpha \cdot \omega > K_X \cdot
  \omega$;

\item If $b_+ = 1$ then the obstruction bundle is of rank zero, and
  hence moduli spaces of symplectic surfaces are \emph{unobstructed}.
\end{itemize}

\noindent The first statement above we knew already, from the duality
$\mathcal{I}(\alpha) = \pm \mathcal{I}(\kappa-\alpha)$.  The second
statement should be interpreted as follows.  View the
invariant $\mathcal{I}(\alpha)$ as a homology class in
$H_{\mathrm{virdim}}(\Gamma(\sss^2; X_r(f)))$, before cutting
  down dimensions; then this invariant is realised by the fundamental
  class of the space of $J$-holomorphic curves for \emph{any} $J \in
  \mathcal{J}$.  Now return to the main theme:

\begin{Thm} \label{nonzero}
Let $X$ be a symplectic four-manifold with $b_+ = 1$ and $b_1 = 0$.
Fix $\alpha \in 
H_2 (X;\zz)$ satisfying $\alpha^2 > K_X \cdot \alpha$ and
$\alpha \cdot \omega > 0$.  Then $\alpha$ contains embedded symplectic
surfaces in $X$.
\end{Thm}

\begin{Pf}
It is enough to prove that for a high degree pencil $f$ on $X$, the
invariant $\mathcal{I}_{(X,f)}(\alpha) \neq 0$.  After twisting by all
the exceptional curves, there is certainly a homotopy class of
sections which yields surfaces in the class $\alpha + \sum E_i$; just
choose a family of $\cdbar$-operators on the restrictions of the line
bundle with this first Chern class to the fibres of $f$, and observe
that an associated projective bundle down the fibres embeds in some $X_r(f)$.
By an earlier result, this homotopy class of sections is unique.

\vspace{0.2cm}

\noindent Fix an almost complex structure $j$ on $P_r(f)$ for which the moduli
space of sections is just a point (a zero-dimensional torus).  Extend
to $J$ on $X_r(f)$ for which the Abel-Jacobi map is holomorphic.  Then
all holomorphic sections of $X_r(f)$ lie over this unique section of
$P_r(f)$, and hence the moduli space is just the space of sections of
some projective bundle.  The conditions on $\alpha$ show that this
projective bundle does indeed have sections (\ref{indextheorem}).  In
the usual way,  this shows that the (compact) space
of sections is a projective space.  Each point condition
$\mathcal{D}(z_i)$ defines a hyperplane in this space, and so the
standard surface count is $\pm 1$.  The result follows.
\end{Pf}

\noindent This is just an instance of the ``wall-crossing formula''.
In general, although the symmetry $\mathcal{I}(\alpha) = \pm
\mathcal{I}(\kappa-\alpha)$ breaks down when $b_+ = 1$, the difference
$|\mathcal{I}(\alpha) - \mathcal{I}(\kappa-\alpha)|$ is given by the
Euler class of an obstruction bundle over $\rot^{b_1(X)/2}$.  Given
(\ref{torusofsections}), one can
prove this using the techniques of this paper, thereby removing the
hypothesis $b_1 = 0$ from the last two applications in (\ref{applications}).


\section{Applications and refinements}

In this section, we shall give the proofs of the applications listed
in the Introduction.  These are due to other authors, and are included
for completeness.  We also give a weak result on homotopy $K3$
surfaces.  At the end of the section, we shall sketch how to
improve the linear constraint $b_+ > 1+b_1$ to $b_+ > 2$, which mildly
improves the scope of the main theorems of this paper and
\cite{SKDIS}.

\subsection{Symplectic consequences}

The first result is an immediate consequence of the main result of
\cite{SKDIS}, and is originally due to Taubes \cite{Taubes}.   

\begin{Cor}
Let $X$ be a symplectic four-manifold with $b_+ > 1+b_1$.  If $X$ is
minimal then $2e(X) + 3\sigma(X) \geq 0$.
\end{Cor}

\begin{Pf}
We may represent $K_X$ by an embedded symplectic surface; by
adjunction, the only components of negative square are $-1$-spheres,
which we exclude by assumption.  So $c_1^2 (X) \geq 0$, and this is
just the inequality given.
\end{Pf}

\noindent If $b_1 (X) = 0$, then the above also shows
that $K_X \cdot [\omega_X] < 0 \Rightarrow b_+ = 1$.
If in fact $b_2 (X) = 1$ then every symplectic form is
rational. Hence, even without proving that one obtains holomorphic
curves from the $\mathcal{I}$-invariants and not just symplectic
surfaces, we recover the theorem of Taubes \cite{Taubes:SWtoGr}:

\begin{Cor}
If $X$ is a symplectic homology projective plane with $K_X \cdot
[\omega] < 0$ then $(X, \omega) \cong (\cc \pp^2, \mu
\omega_{FS})$ for some $\mu > 0$.
\end{Cor}

\begin{Pf}
According to (\ref{nonzero}) the generator for the homology contains
smooth symplectic surfaces.  These must be connected, by irreducibility of
the homology class, and then the adjunction formula shows that such a
surface must in fact be a sphere.  Since the virtual dimension of the
moduli space was $2$, there is in fact such a sphere through any two
generic points of $X$.  Then, as in Gromov \cite{Gromov}, one uses
such a family of
spheres to construct an explicit diffeomorphism from $X$ to $\cc
\pp^2$ which pulls back the standard symplectic form.
\end{Pf}

\noindent This was the original motivation for the section.  We can
extend that result as follows: the following is due to Ohta and Ono
\cite{Ohta-Ono}. 

\begin{Cor}
Let $X$ be a symplectic four-manifold with $b_1 = 0$ and with $K_X
\cdot [\omega_X] < 0$.  Then $X$ is diffeomorphic
to $\sss^2 \times \sss^2$ or to $\cc \pp^2 \sharp n
\overline{\cc\pp}^2$ with $n \leq 8$.
\end{Cor}

\begin{Pf}
Again it is enough to construct a symplectic sphere of non-negative square.
Then work of McDuff provides a diffeomorphism to one of the listed del
Pezzo surfaces (and in fact work of Lalonde-McDuff 
on the classification of symplectic forms on rational manifolds pins down
the precise symplectic structure).  Note that 
their work does not rely on gauge theory!  Now
if $K_X \cdot \omega < 0$ then for any $\alpha$
with $\alpha \cdot \omega > 0$ we have that $\iota(\alpha) \cdot
[\Fibre] > 2g-2$; this is a trivial consequence of adjunction.  It
follows, by the remarks before Theorem (\ref{nonzero}),
that if $\alpha$ satisfies $\alpha^2 > K\cdot \alpha$ then
$\alpha$ will have embedded symplectic representatives through any
point.  At least one
component of each will be a sphere provided $\alpha^2 + K_X \cdot
\alpha < 0$, and this underlies the result.

\vspace{0.2cm}

\noindent Since $K_X = -\lambda [\omega]$ we know $c_1^2 (X) > 0$.
However, since $b_+ = 1$ and $b_1 = 0$, then 

$$2e(X) + 3\sigma(X) \ = \ 4+2(1+b_-) + 3(1-b_-) = 9-b_-.$$

\noindent It follows that $b_- \leq 8$.  The Hasse-Minkowski
classification of intersection forms shows that the only possible intersection
forms are $\zz \langle 1 \rangle \oplus b_- \zz \langle -1 \rangle$,
in the odd case, and the rank two hyperbolic form in the even case.
In both
situations, the generator of the positive summand $\alpha =
(1,\underline{0})$ is quickly checked to satisfy all the required
conditions, and this homology class contains symplectic spheres. 
\end{Pf}

\noindent Finally, we note that these methods give rise to
constructions of symplectic 
forms on four-manifolds with $b_+ = 1$.  There are two standard ways
of building symplectic forms on four-manifolds:

\begin{itemize}
\item By \emph{integration}; find a large family $\mathcal{S}$ of
  irreducible homologous surfaces which
  cover $X$ and all have pairwise locally positive strictly positive
  intersections.  Define a current by taking a two-form $\chi$ to the
  number $\int_{u \in \mathcal{S}} \int_{S_u} \chi$; under fairly
  general circumstances, this defines a symplectic form (cf. Gompf's
  work in \cite{Gompf} and \cite{GompfGokova}) in the class Poincar\'e
  dual to the homology class of any surface.

\item By \emph{inflation};  be given one symplectic form $\omega$ on
  $X$, and a connected symplectic surface $S \subset X$, and then
  deform $\omega$ by adding positive forms supported in the normal
  bundle of $S$ (cf. McDuff's \cite{McDuff:DeftoIso} and Biran's
  \cite{Biran}).  One obtains forms in the class $[\omega] + tPD[S]$
  for all $t \geq 0$.
\end{itemize}

\noindent The second method is just a special case of the first: the
normal bundle of the surface $S$ is filled by homologous symplectic
surfaces.  We can use either method here; the latter is mildly
simpler.  This result, and generalisations which also follow from our
work, was noted by Li and Liu \cite{Li-Liu}.

\begin{Cor}
Suppose that $X$ is minimal and $b_+ = 1$, $b_1 = 0$.  If
$K^2 > 0$ and $K \cdot \omega > 0$ then the canonical class contains
symplectic forms.
\end{Cor}

\begin{Pf} For large enough $N$, and
arguments as above, we find that the homology class $NK-[\omega]$
contains symplectic surfaces in $X$.  An easy argument with the
intersection form (and the minimality assumption) ensures that
these surfaces are connected.  Then inflate!  
\end{Pf}

\noindent  This finishes our
treatment of the applications (\ref{applications}); but we will end by
discussing an application that we cannot yet complete.
In the proof of (\ref{Serreduality}), we used the
Brill-Noether theory to get control on the 
geometry of the map $X_r(f) \rightarrow P_r(f)$.  For $r$ neither very
large nor small
compared to $g$, say $r=g-1$, we have no such control.  Hence we cheat
and make a definition.  Every Riemann surface of genus $g$ has an
associated $\Theta$-divisor, the image of $\Sym^{g-1}(C) \rightarrow
\Pic_{g-1}(C)$.  The image is some subvariety of the Picard torus, and
for a Zariski open set $\mathcal{U}$ in moduli space the deformation
type of the 
pair $(\Pic_{g-1}, \Sym^{g-1})$ will vary locally constantly.  It is
easy to check that the complement $\mathcal{Q}$ of this open set has
divisorial components. 

\begin{Defn}
A symplectic four-manifold $X$ is $\Theta$-positive if for all 
Lefschetz pencils 
$f$ of high enough degree, the associated sphere in moduli space
$\mgbar$ meets $\mathcal{Q}$ locally positively. 
\end{Defn}

\noindent Algebraic positivity is proved in
\cite{ivanmodulidivisor}, but our techniques do not yet permit ``Whitney
moves'' in $\mgbar$ to cancel excess intersections. At least for the
class of $\Theta$-positive manifolds, we can reprove a nice result of
Morgan and Szabo, which follows from \cite{Morg-Sz}.  We shall just
sketch the idea of the proof.

\begin{Prop}
Suppose $X$ is a $\Theta$-positive, simply-connected symplectic
four-manifold with $c_1(X) = 0$. Then $X$ is a homotopy $K3$ surface.
\end{Prop}

\begin{Pf}[Sketch]
It would be enough to prove the following (which is Morgan
and Szabo's theorem):  for $X$  spin symplectic with $b_1 = 0$ and
$c_1^2 = 0$ then $Gr_X (K_X / 2)$ is odd iff $b_+ = 3$.  This is the
mod two version of the computation given at the end of section two.
(There are examples which show that in the
symplectic category,
the Gromov invariant of $K_X / 2$ can be any even or odd number for each
fixed homotopy type, depending on whether $b_+ > 3$ or $b_+ = 3$.)
Let us re-cast the last stage of that earlier argument.  The map

$$H^0 (K_X / 2) \times H^0 (K_X / 2) \ \rightarrow \ H^0 (K_X)$$

\noindent induces a map on projective spaces $\pp^a \times \pp^a
\rightarrow \pp^{2a}$ which can be described as follows:  $\pp^a =
\Sym^a (\pp^1)$ and the map is just addition of divisors.  The
K\"ahler surfaces of the given homotopy type are elliptic; each
(half-)canonical divisor is a collection of fibres, determined by a
set of points on the base 
$\pp^1$ of the elliptic fibration.  Now the Gromov invariant of $K_X /
2$ is just the \emph{degree} of the map $+: \pp^a \times \pp^a
\rightarrow \pp^{2a}$, which is even provided the image has positive
dimension, i.e. if $b_+ > 3$. 

\vspace{0.2cm}

\noindent Be given a spin symplectic manifold.  If  $c_1^2 =
0$ we know that $b_+ \geq 3$, and if also $b_1 = 0$ then the index of
the $\cdbar$-operator on the Picard bundle is necessarily negative.
Split the exceptional sections into two collections of equal numbers
of spheres $E_{\mathcal{A}}$ and $E_{\mathcal{B}}$.  We obtain two
sections of $P_{g-1}(f)$, corresponding to $K_X/2 + E_{\mathcal{A}}$
and $K_{X'} - (K_X/2 + E_{\mathcal{A}}) \ = \ K_X/2 +
E_{\mathcal{B}}$, and by generically extending some $J_{\Pic}$ away from
these we can assume these are the only holomorphic sections of the
Picard bundle.  If $X$ is $\Theta$-positive then we can find a $J$ on
$X_{g-1}(f)$ which maps holomorphically to the Picard bundle with
structure $J_{\Pic}$.  Then the moduli space of holomorphic
sections of 
$X_{g-1}(f)$ in the class $s_{\mathcal{A}}$ is a projective space
$\pp^a$ and in the class $s_{\mathcal{B}}$ it's a projective space
$\pp^a$ of the same dimension.  Moreover, we have a map

$$+: \pp^a \times \pp^a \ \rightarrow \ \pp^N; \qquad N \ = \
[b_+ - 3]/2.$$ 
  
\noindent This multiplies holomorphic curves in the obvious way to
give a section of the canonical bundle on the fibre: the projective
space on the right is just the moduli space of holomorphic sections of
$X_{2g-2}(f)$ for an almost complex structure standard near the
zero-sections.  On each fibre,
then, the $+$ map is given
by the natural map $\Sym^{g-1} \times \Sym^{g-1}
\rightarrow \Sym^{2g-2}$.  But, if $a>0$, this last map has \emph{even
  degree} onto its image;
it factors through $\Sym^2 (\pp^a)$ for instance.  The standard
surface count for the class $K_X/2$ can be obtained from the
obstruction bundle over $\pp^a$ which can be described in terms of the
obstruction bundle over its image in $\pp^N$ via pullback by the
map $+$.  The 
evenness of the degree of this map translates to saying that the Euler
class is even, and we finally obtain that when $b_+ > 3$ the standard
surface count for $K_X / 2$ is even.  But this, together with
\cite{SKDIS}, implies that $K_X \neq 0$; hence $\Theta$-positivity implies
the Morgan-Szabo theorem.  
\end{Pf}

\noindent It would be interesting to see if one could
use the geometry of the theta-divisor, or of harmonic maps to the
moduli space of curves, to prove the positivity always holds.


\subsection{A better linear constraint}

Throughout the last sections, we have worked with complex
structures on bundles and Jacobians generic away from certain specified
loci.  The factor $b_1 (X)$ in all the various inequalities
occurs because (\ref{torusofsections} aside) we have ignored the
topological structure -  in the form
of trivial subrepresentations of the homological monodromy - coming from
line bundles on the four-manifold.  We shall now indicate one naive route to
taking account of the extra structure; in doing so we will
improve the linear
constraint in the main theorem of \cite{SKDIS} from
$b_+ > 1+b_1$ to the slightly improved constraint $b_+ > 2$.  A similar
analysis would  yield the same improvement for the main theorem of this
paper.  This
still falls short of Taubes, however, so the discussion should be
regarded as somewhat parenthetical.  We include it only to indicate
that the problems arising here are related to the failure of the Hard
Lefschetz theorem for symplectic manifolds.  Indeed, if $b_1 (X)$ is odd, then
we can't ``share out'' the homology equally between $H^{1,0}(\Sigma)$
and $H^{0,1}(\Sigma)$ for a hyperplane $\Sigma$.  But this introduces
an asymmetry into the construction: the projectivisation of the first vector
space is the
fibre of the bundle of fibrewise canonical forms, and the second
vector space is the fibre of the tangent bundle to the Picard
fibration at the zero-section.  Presumably a different adaptation will
avoid this hiccup, but it sheds some light on the different role of
the first homology of the four-manifold in our treatment and that of Taubes.
The argument is very closely related to (\ref{torusofsections}).

\vspace{0.2cm}

\noindent As 
piece of notation, write $R$ for $b_1 (X) / 2$ when $b_1$ is even and
for $(b_1 (X) - 1)/2$ when $b_1$ is odd.  We adopt the notations of
\cite{SKDIS} wherever not already defined.  In particular, $f_* K$
denotes the vector bundle of fibrewise canonical forms.  That is, if the
bundle $W$ has fibre canonically $H^0 (\Sigma_t; K_{\Sigma_t})$ over $t
\in \pp^1$, then $f_* K = W \otimes K_{\pp^1}$.  The following
improves (\ref{indextheorem}).

\begin{Prop} \label{betterindex}
If $b_1(X)$ is even then there is a holomorphic structure on $f_* K$
for which the space of holomorphic sections has dimension $(b_+ - 1)/2$.
If $b_1 (X)$ is odd there is a holomorphic structure for which the space
of holomorphic sections has dimension $(b_+-2) / 2$.
\end{Prop}

\begin{Pf}
For a K\"ahler surface there is a  subbundle of $f_* \omega_{X/B}$ with fibre
at a point $t$ given by $H^{1,0}(X) = H^0 (X, \Omega^1_X)$.  The
subbundle arises from the obvious restriction map on holomorphic
differentials. 
Moreover there is a decomposition $H^1 (X, \cc) = H^{1,0} \oplus
H^{0,1}$.  The vector space on the LHS carries a trivial monodromy
representation for any (Lefschetz) fibration, and hence we have a 
holomorphically trivial subbundle of $f_* \omega_{X/B}$ over
$\pp^1$.  The rank of the subbundle is $R$.  The proposition relies on
the analogue of this for a general symplectic pencil.

\vspace{0.2cm}

\noindent Using a family of metrics on the fibres
we can write, at any point $t \in \sss^2$, the cohomology group

$$H^1 (\Sigma_t, \cc) \ = \ H^{1,0} (\Sigma_t) \oplus H^{0,1}
(\Sigma_t) \ = \ H^0 (\Sigma_t, K_{\Sigma_t}) \oplus H^0 (\Sigma_t,
K_{\Sigma_t})^*.$$

\noindent We have an inclusion $H^1 (X, \cc) \subset H^1
(\Sigma_t, \cc)$ as the subgroup of monodromy invariants, but now
there may be no non-trivial intersection of $H^1 (X, \cc)$ with
$H^{1,0}(\Sigma_t)$.  We can change the choice of almost complex
structures on the fibres to avoid this.  Recall that a complex
structure on the real vector space $H^1 (\Sigma, \rr)$ is precisely a
decomposition of the complexification $H^1 (\Sigma, \cc) = H^1
(\Sigma, \rr) \otimes \cc$ into two half-dimensional subspaces which
are conjugate under the action of the conjugation on the tensor
product factor.  Given one reference complex structure $H^1 (\Sigma,
\cc) = E \oplus \overline{E}$, for instance a splitting into
holomorphic and antiholomorphic forms, any other is
given by a complex linear homomorphism $E \rightarrow \overline{E}$
whose graph defines one of the two new decomposing subspaces.  Now the
subspace $H^1 (X, \cc) \subset H^1 (\Sigma, \cc)$ is preserved by the
complex conjugation, though no non-trivial complex subspace is fixed
pointwise.  It follows that we can choose a homomorphism in $Hom (E,
\overline{E})$ for which an $R$-dimensional subspace of $H^1 (X, \cc)$
lies inside the new space of holomorphic forms.  Equivalently, we can
choose a fibrewise family of complex structures on the fibres of a
Lefschetz pencil to ensure that the intersection of the trivial
complex bundle with fibre $H^1 (X, \cc)$ and the bundle $f_* \omega_{X
  / \sss^2}$ has rank $R$.

\vspace{0.2cm}

\noindent It follows that there is a sequence of topological bundles

\begin{Eqn} \label{trivialinside}
0 \rightarrow \underline{\cc}^R \rightarrow f_* \omega_{X/\sss^2}
\rightarrow Q \rightarrow 0
\end{Eqn}

\noindent where $Q$ is defined as the cokernel: this gives a sequence

$$0 \rightarrow (K_{\pp^1})^{\otimes R} \rightarrow f_*K \rightarrow
Q' \rightarrow 0.$$

\noindent Suppose we choose a complex structure on the vector bundle
$f_* K$ which makes the subbundle of the sequence holomorphic and
which extends that connexion generically to $Q'$.  Then
from the long exact sequence in cohomology we find that $H^0 (f_* K) =
\ind_Q (\cdbar)$ which is easily computed to be as claimed in the
proposition. 
\end{Pf}

\noindent Given this, one can find holomorphic sections of the
projective bundle whenever $b_+ > 2$.  To run the rest of the
computation of $\mathcal{I}(\kappa)$ requires one or two further
modifications.  The principal of these is a new definition for the
complex structure $J_{\mathrm{ext}}$ compatible with the zero-sections
(or more generally compatible with duality).  We start with the linear
complex structure on $\pp (f_* K)$ provided by (\ref{betterindex}).
Inside the fibration of degree zero Jacobian varieties we have a
subfibration $\Pic_f ^X$ of tori of dimension $R$, using the sequence
(\ref{trivialinside}) and the duality given at the end of Remark
(\ref{exactseqs}).   These
fibres are Picard varieties for line bundles on the symplectic
manifold $X$.  There
is an analogous subfibration inside all the higher degree Picard
fibrations, defined up to translation.  We can
choose a complex structure on $\Sym^{(2g-2)}_f$ to make not only
$\pp(f_* K)$ a holomorphic subset but also the preimage of this entire
subfibration of $P_{2g-2}(f)$.  Choose a generic such almost
complex structure.

\begin{Lem}
Let $X$ be a symplectic four-manifold with $2 < b_+ (X) \leq 1+b_1 (X)$. 
For the $J_{\mathrm{ext}}$ described above, the whole moduli space of
sections of $\Sym^{(2g-2)}_f$ is still the projective space of
sections of $\pp (f_* K)$.
\end{Lem}

\begin{Pf}
Suppose as before we have two sections $u_1, u_2$ of the fibration of
symmetric products.  As in (\ref{torusofsections}), the difference
$\tau(u_1) (\tau(u_2))^{-1}$  is well-defined 
as a section of the bundle $R^1 f_* \mathcal{O} \rightarrow \sss^2$.
The sections of this are dual to the constant sections of the trivial
subbundle of $f_* \omega$.  It follows that although not all sections
of $\Pic_f$ co-incide with the image under $\tau$ of the projective
bundle $\pp(f_* K)$, all sections are constant translates of this
section still lying inside the subfibration $\Pic_f ^X$.
Hence all the holomorphic sections of $\Sym^{(2g-2)}_f$ are in fact
holomorphic sections of $\pp (f_* K)$ or of a projective bundle of
rank $g-2$ over $\sss^2$.  But the index of the $\cdbar$-operator on
all of these other projective bundles is still negative, by the
assumption on the Betti numbers of $X$.  Hence for a generic extension
of the almost complex structure from $\pp (f_* K)$ to the rest of the
fibration $\tau^{-1} (\Pic_f ^X)$, none of these other projective
bundles have any sections; not all sections of $P_r(f)$ lift to the
fibration of symmetric products.  The result follows.
\end{Pf}

\noindent Note that it is easier to ``fill in'' the remaining
four-manifolds with $2< b_+ \leq 1+b_1$ than to give a treatment for
all four-manifolds with $b_+ > 1$ in one step.  There is one last
alteration required in the proof.
In the long exact sequence in cohomology which underlies the
obstruction computation (\ref{obstructionsequence}), the first
$H^1$-term no longer vanishes,
  since there is a non-trivial obstruction bundle for sections when
  viewed as lying inside $\pp (f_* K)$.  However, the preceding map in
  the exact sequence with image this obstruction space is an
  isomorphism, and the
  rest of the argument proceeds much as before.   


\bibliographystyle{alpha}
\bibliography{main}

\end{document}